%% file: paper_siam.tex
\documentclass[final]{siamltex}

\usepackage{latexsym,amsmath,amssymb,amsfonts,amscd,multirow,dsfont}
\usepackage{epsfig,verbatim,epstopdf,graphics}
\usepackage{bm,color}
\usepackage{subfigure}
\usepackage{mathtools}
\usepackage{enumitem}
\usepackage[title]{appendix}

\newtheorem{thm}{Theorem}[section]

\newtheorem{lem}[thm]{Lemma}
\newtheorem{exa}[thm]{Example}
\newtheorem{rem}[thm]{Remark}

\newcommand{\bk}{\mathbf{k}}
\newcommand{\bx}{\mathbf{x}}

\newfont{\iams}{msbm9}

\newcommand{\commentbis}[1]{}
\newcommand{\be}{\begin{eqnarray}}
\newcommand{\ee}{\end{eqnarray}}
\newcommand{\beno}{\begin{eqnarray*}}
	\newcommand{\eeno}{\end{eqnarray*}}

\newcommand{\barr}[1]{\begin{array}{#1}}
	\newcommand{\earr}{\end{array}}

\newcommand{\beq}{\begin{equation}}
\newcommand{\eeq}{\end{equation}}
\newcommand{\beqa}{\begin{eqnarray}}
\newcommand{\eeqa}{\end{eqnarray}}

\newcommand{\bV}{{\bf V}}
\newcommand{\bn}{{\bf n}}

\newcommand{\bP}{{\bf P}}
\newcommand{\bQ}{{\bf Q}}
\newcommand{\bR}{{\bf R}}

\newcommand{\bzero}{\mathbf{0}}
\newcommand{\bone}{\mathbf{1}}
\newcommand{\bl}{\mathbf{l}}
\newcommand{\bi}{\mathbf{i}}

\newcommand{\bj}{\mathbf{j}}
\newcommand{\bW}{\mathbf{W}}

\newcommand{\bK}{\mathbf{K}}
\newcommand{\ba}{{\bm{\alpha}}}
\newcommand{\bb}{{\bm{\beta}}}
\newcommand{\bq}{\mathbf{q}}
\newcommand{\vb}{v_{\mathbf{i},\mathbf{l}}^\mathbf{j}}
\newcommand{\vbp}{v_{\mathbf{i}',\mathbf{l}'}^{\mathbf{j}'}}
\newcommand{\sgv}{\widehat{\mathbf{V}}_N^k}

\newcommand{\mwk}{\mathcal{W}^{\infty,k+1}}
\newcommand{\mQ}{\mathcal{Q}_{\bl,\bl'}^{\bj,\bj'}}
\makeatletter
\newcommand{\Rmnum}[1]{\expandafter\@slowromancap\romannumeral #1@}
\makeatother
\newcommand{\vertiii}[1]{{\left\vert\kern-0.25ex\left\vert\kern-0.25ex\left\vert #1 
		\right\vert\kern-0.25ex\right\vert\kern-0.25ex\right\vert}}

\allowdisplaybreaks[1]

\graphicspath{{./}{./figures/}}

\title
{A sparse grid discontinuous Galerkin method for the high-dimensional Helmholtz equation with variable coefficients}

\author{
	Wei Guo
	\thanks{Department of Mathematics and Statistics, Texas Tech University,
		Lubbock, TX, 70409. Research is supported by NSF grants DMS-1620047, DMS-1830838.
		{\tt weimath.guo@ttu.edu}}
}

\begin{document}
\maketitle

\begin{abstract}
	%
	The simulation of high-dimensional problems with manageable computational resource represents a long standing challenge.  In a series of our recent work \cite{wang2016sparse,guo2016sparse1,guo2017adaptive,tao2018sparseguo}, a class of sparse grid DG methods has been formulated for solving various types of partial differential equations in high dimensions. By making use of the multiwavelet tensor-product bases on sparse grids in conjunction with the standard DG weak formulation, such a novel method is able to significantly reduce the computation and storage cost compared with full grid DG counterpart, while not compromising accuracy much for sufficiently smooth solutions. In this paper, we consider the high-dimensional  Helmholtz equation with variable coefficients and demonstrate that for such a problem the efficiency of the sparse grid DG method can be further enhanced by exploring a semi-orthogonality property associated with the multiwavelet bases, motivated by the work \cite{pflaum1998multilevel,pflaum2016sparse,hartmann2018prewavelet}. 
	The detailed convergence analysis shows that the modified sparse grid DG method attains the same order accuracy, but the resulting stiffness matrix is much sparser than that by the original  method, leading to extra computational savings. Numerical tests up to six dimensions are provided to verify the analysis. 

\end{abstract}


\section{Introduction}


Discontinuous Galerkin (DG) methods that make use of discontinuous functions as approximations have been extensively studied for solving partial differential equations (PDEs) over the last few decades and become rather mature for various applications \cite{cockburn2000development}. It is generally understood that DG methods  are more flexible compared with continuous finite element methods  due to the lack of continuity requirement, and thus enjoying many attractive properties. However, the situation changes when the dimension of the underlying problem becomes large. Traditional grid-based methods including DG methods suffer from the \emph{curse of dimensionality} \cite{bellman1961adaptive}, which describes the scenario that the complexity and memory storage of an algorithm grows exponentially with the dimension of the underly problem for a given level of accuracy. Moreover, compared 
with other high order schemes, DG methods often require more degrees of freedom (DOF), and hence the effect of the curse of dimensionality is more significant, making such methods uncompetitive for high-dimensional calculations.

Driven by the need for a
method that is able to retain the attractive properties of DG methods and yet requires feasible computational
and storage cost for high-dimensional simulations, based on the sparse grid approach \cite{zenger1991sparse,bungartz2004sparse, garcke2013sparse}, a class of novel DG methods has been developed for solving various types of high-dimensional PDEs in a series of work \cite{wang2016sparse,guo2016sparse1,guo2017adaptive,tao2018sparseguo,tao2018sparse,liu2019krylov}. For such a method, the underlying DG solution is represented by a set of orthonormal tensor-product multiwavelet bases \cite{alpert1993class,wang2016sparse,guo2016sparse1}. Meanwhile,  following the standard sparse grid philosophy, rather than including all the anisotropic tensor-product bases, only the ones with significant contribution to the approximation accuracy are chosen, leading to immense reduction in cost complexity when the dimension $d$ is large. In particular, the number of DOF scales as $\mathcal{O}(h^{-1}| \log h|^{d-1})$ instead of $\mathcal{O}(h^{-d})$, where $h$ denotes the mesh size in each direction; while the approximation property of the proposed sparse grid DG method is proven only slightly deteriorated for smooth solutions through both theoretical and numerical verifications, see \cite{wang2016sparse,guo2016sparse1}. 


In this paper, we are concerned with the Helmholtz equation with variable coefficients, which arises in many applications in science and engineering such as the representation of the solution of Schr\"{o}dinger equation \cite{pflaum2016sparse}. In \cite{wang2016sparse}, we developed an efficient high order sparse grid DG method based on the  interior penalty DG (IPDG) weak formulation \cite{baker1977finite, douglas1976interior, wheeler1978elliptic, arnold1982interior,arnold2002unified} for solving general elliptic equations, and such a method can be directly applied to simulate the problem. Due to the reduction of DOF, the resulting algebraic linear system enjoys a much smaller dimension than that by the traditional DG method, leading to computational savings in high dimensions. However,  the linear system also becomes denser especially in the variable-coefficient case because of the hierarchical nature of the multiwavelet bases, and thus impeding the efficiency advantage of the sparse grid approach. To address the issue, we would like to explore the semi-orthogonality property
associated with the orthonormal multiwavelet bases, which can be thought of as a special type of orthogonality of the bases related to the underlying bilinear formulation.
The idea of making use of semi-orthogonality on sparse grids was proposed in \cite{pflaum1998multilevel,pflaum2016sparse,hartmann2018prewavelet} in the  continuous Ritz-Galerkin discretization framework. In particular, by using the prewavelet as the bases, the standard bilinear formulation is carefully modified according to the semi-orthogonality property, aiming to sparsify the stiffness matrix for the variable-coefficient equations. 
The analysis in \cite{pflaum2016sparse} further shows that such a modification preserves the accuracy of the original method under some extra smoothness requirement of the coefficient. In this work, we demonstrate that the modification of the bilinear form based on semi-orthogonality is also  effective for our sparse grid IPDG method with multiwavelet bases. In the theoretical aspect, we refine the analysis in \cite{pflaum2016sparse} and show that the error incurred by modifying the bilinear form is one order higher than the projection error for sufficiently smooth problems, making the modified method produce virtually the same numerical result as the original sparse grid IPDG method; while the associated linear system is much sparser,  offering extra computational and storage savings. Furthermore, under the sparse grid DG framework, the proposed methodology has the potential to be extended to many other variable-coefficient PDEs.

The rest of this paper is organized as follows.
In Section \ref{sec:IPDG},  we review the fundamentals of the IPDG on sparse grids for solving the Helmholtz equation developed in \cite{wang2016sparse,guo2016sparse1}.  In Section \ref{sec:mIPDG}, we formulate the modified scheme based on the semi-orthogonality property and perform an error analysis. In Section \ref{sec: SGDG numerical}, numerical results in multi-dimensions (up to $d=6$) are provided to validate the accuracy and performance of the  method. Conclusions  and future work are given in Section \ref{sec:Con}.

\section{IPDG Method on Sparse Grids for the Helmholtz Equation}
\label{sec:IPDG}

In this section, we review the construction together with several key theoretical results of the  IPDG method on sparse grids for solving the following Helmholtz equation,
\begin{align}
-\Delta u + c u =& f\quad \text{in}\quad \Omega=[0,1]^d,\label{eq:helmd}\\
u=&0 \quad\text{on} \quad \partial \Omega,\notag
\end{align}
where $c$ is a non-negative function and the domain is the $d$-dimensional unit box.

\subsection{DG Finite Element Spaces on Sparse Grids}
\label{sec:SGapproximation}
In this subsection, we briefly review the recent development on the DG finite element approximation space on sparse grids in \cite{wang2016sparse,guo2016sparse1}, which plays a key role in dimension reduction for the proposed DG method.


%


We first introduce the hierarchical
decomposition of DG approximation space in one dimension.
%
%
Without loss of generality, we consider the domain $\Omega=[0,1]$ and define the $l$-th level grid $\Omega_l$, composed of non-overlapping intervals $I_{l}^j=(2^{-l}j, 2^{-l}(j+1)]$, $j=0, \ldots, 2^l-1$ with uniform cell size $h_l=2^{-l}.$ 
%
%
Denote by
$$V_l^k=\{v: v \in P^k(I_{l}^j),\, \forall \,j=0, \ldots, 2^l-1\}$$
the collection of piecewise polynomials of degree at most $k$ defined on grid $\Omega_l$.  Note  that there exists a nested structure with respect to different values of $l$:
$$V_0^k \subset V_1^k \subset V_2^k \subset V_3^k \subset  \cdots,$$
which allows us to define the increment space $W_l^k$, $l=1, 2, \ldots $ as the orthogonal complement of $V_{l-1}^k$ in $V_{l}^k$ with respect to the standard $L^2$ inner product on $\Omega$, i.e.
\begin{equation*}
V_{l-1}^k \oplus W_l^k=V_{l}^k, \quad W_l^k \perp V_{l-1}^k.
\end{equation*}
By letting  $W_0^k:=V_0^k$, which consists of all polynomials of up to degree $k$ on the entire domain $[0,1]$,
we arrive at the hierarchical decomposition of  the standard piecewise polynomial space $V_N^k$, $N\ge1$ on grid $\Omega_N$ as \begin{equation}
\label{eq:hr}
V_N^k=\bigoplus_{0 \leq l\leq N} W_l^k.
\end{equation}



In light of \eqref{eq:hr}, we are able to represent a function in $V_N^k$ by making use of the orthonormal multiwavelet bases of space $W_l^k$ developed in \cite{alpert1993class}.  We denote such orthonormal bases  of $W_l^k$ as
$$
v_{i, l}^{j}(x), \quad i=1, \ldots, k+1, \quad j=0, \ldots, \max(0,2^{l-1}-1).
$$ 
The details about the construction of the multiwavelet bases can be found in \cite{alpert1993class,wang2016sparse}.

It is worth noting that space $W_l^k$ can be constructed by means of projection operators as well, which will be helpful in the analysis below. Denote $P_l^k$ as the standard $L^2$  projection operator from  $L^2(0,1)$ onto $V_l^k$, i.e., $V_l^k=P_l^k(L^2(0,1))$.
Then, we define the
increment projector,
\[Q_l^k: =
\begin{cases}
P_l^k - P_{l-1}^k,  & \textrm{if} \quad l \ge1,\\
P_0^k, & \textrm{if} \quad l = 0.
\end{cases}
\]
Evidently,
$$W_l^k= Q_l^k( L^2(0,1)).$$
To summarize, we have obtained the hierarchical decomposition of the piecewise polynomial space $V_N^k$ as well as the corresponding projection operator $P_N^k$:
$$V_N^k=\bigoplus_{\substack{0 \leq l\leq N}} W_l^k,  \ P_N^k=\sum_{\substack{0 \leq l\leq N}} Q_l^k.$$

Then, we review the construction of approximation spaces in multi-dimensions. Consider the domain $\bx=(x_1, \ldots, x_d) \in \Omega=[0,1]^d$. We first recall some notation on the norms and operations of multi-indices in $\mathbb{N}_0^d$, where $\mathbb{N}_0$ denotes the set of nonnegative integers. For $\ba =(\alpha_1, \ldots, \alpha_d) \in \mathbb{N}_0^d,$ we define the $l^1$ and $l^\infty$ norms
$$
|\ba|_1:=\sum_{m=1}^d \alpha_m, \qquad   |\ba|_\infty:=\max_{1\leq m \leq d} \alpha_m,
$$
the component-wise arithmetic operations
$$
\ba \pm \bb :=(\alpha_1 \pm \beta_1, \ldots, \alpha_d \pm\beta_d),  \qquad 2^\ba:=(2^{\alpha_1}, \ldots, 2^{\alpha_d}),
$$
the relational operator
$$
\ba \leq \bb \Leftrightarrow \alpha_m \leq \beta_m, \, \forall m,
$$
and
$$
\max(\ba,\bb) := (\max(\alpha_1,\beta_1),\ldots,\max(\alpha_d,\beta_d)),\quad \min(\ba,\bb) := (\min(\alpha_1,\beta_1),\ldots,\min(\alpha_d,\beta_d)).
$$
Denote $\bzero:=(0,\ldots,0)\in\mathbb{N}_0^d $ and $\bone:=(1,\ldots,1)\in\mathbb{N}_0^d$. 

Given a multi-index $\bl=(l_1, \ldots, l_d) \in \mathbb{N}_0^d$ indicating the levels of the mesh in a multivariate sense, we define the grid $\Omega_\bl$ consisting of
non-overlapping elementary cells  $I_\bl^\bj=\{\bx: x_m \in (2^{-l_m}j_m, 2^{-l_m}(j_m+1)),\,m=1,\ldots,d\}$, $\bzero \le \bj \le 2^{\bl}-\bone$.
%
%
On grid $\Omega_\bl$, we define the tensor-product increment space
%
$$\bW_\bl^k:=W_{l_1, x_1}^k \times W_{l_2, x_2}^k \cdots \times W_{l_d, x_d}^k,$$
where $W_{l_m, x_m}^k$ corresponds to the one-dimensional  (1D) increment space in the $m$-th direction.
Based on the 1D hierarchical decomposition \eqref{eq:hr}, we yield
$$\bV_N^k:=V_{N,x_1}^k \times V_{N,x_2}^k \cdots \times V_{N,x_d}^k=\bigoplus_{|\bl|_\infty \leq N} \bW_\bl^k,$$
where $\bV_N^k$ is indeed the traditional tensor-product DG approximation space on grid $\Omega_N:=\Omega_{(N,\ldots,N)}$, i.e., having 2$^N$ cells in each dimension. $\bV_N^k$ will be called the full grid space hereafter.

Define the index set $\Theta_{\bl}=\{\bj\in \mathbb{N}_0:\bzero \le \bj \le\max(\bzero, 2^{\bl-\bone}-\bone)\}$. The basis functions for $\bW_\bl^k$ are constructed by a tensor product of the 1D multiwavelet bases
\begin{equation}\label{eq:basisd}
v_{\bi, \bl}^{\bj}(\bx) := \prod_{m=1}^d {v_{i_m, l_m}^{j_m}(x_m)}, \quad \bj\in \Theta_{\bl},\, \mathbf{1}\le\bi\le\bk+\mathbf{1},
\end{equation}
where $\bk:=(k,\ldots,k)\in\mathbb{N}^d_0$. 
Note that they form a set of orthonormal bases due to the property of the 1D bases.
The sparse grid DG finite element approximation space is defined as
$$\widehat{\bV}_N^k:=\bigoplus_{\substack{ |\bl|_1 \leq N} }\bW_\bl^k.$$
Evidently, $\widehat{\bV}_N^k$ is a subset of $\bV_N^k$. More importantly, its number of DOF scales as $O((k +
1)^d2^NN^{d-1})$, as opposed to  exponential dependence on $Nd$ for the full grid space $\bV_N^k$, see \cite{wang2016sparse}. This is the key for computational savings in high dimensions.

The standard $L^2$ projection operator onto space $\bW_l^k$ can be naturally written as
$$
\bQ_\bl^k := Q^k_{l_1,x_1}\otimes\cdots\otimes Q^k_{l_d,x_d},
$$
where ${Q}^k_{l_m, x_m}$ denotes the 1D projection operator ${Q}^k_{l_m}$ in the $m$-th dimension. Consequently, we write the $L^2$ projection operator onto the sparse approximation space $\widehat{\bV}_N^k$ as
\begin{equation}
\label{eq:sgproj}
{\widehat{\bP}}_N^k:=\sum_{\substack{ |\bl|_1 \leq N}} \bQ_{\bl}^k.
\end{equation}

\subsection{IPDG Method on Sparse Grids}
\label{sec:sparseIP}




To formulate the IPDG method on sparse grids, we start with  introducing some standard notation about jumps and averages for   piecewise functions defined on grid $\Omega_N$. Let $\Gamma:=\bigcup_{T \in \Omega_N} \partial_T$ be the union of the boundaries for all the elements in $\Omega_N$ and $S(\Gamma):=\Pi_{T\in \Omega_N} L^2(\partial T)$ be the set of $L^2$ functions defined on $\Gamma$. For any $q \in S(\Gamma)$ and $\bq \in [S(\Gamma)]^d$,   their   averages $\{q\}, \{\bq\}$ and jumps $[q], [\bq]$ on the interior edges are defined as follows,
\begin{flalign*}
[ q] \  =\  \, q^- \bm{n}^- \, +  q^+ \bm{n}^+, & \quad \{q\} = \frac{1}{2}( q^- + 	q^+), \\
[ \bq] \  =\  \, \bq^- \cdot \bm{n}^- \, +  \bq^+ \cdot \bm{n}^+, & \quad \{\bq\} = \frac{1}{2}( \bq^- + \bq^+),
\end{flalign*}
where $\bn$ denotes the unit normal, and `$-$' and `+' represent the directions of the underlying vector pointing to interior and exterior at an element edge $e$,  respectively. If $e$ is a boundary edge, then we let
$$[q] \  =q \bm{n}, \quad \{\bq\} =\bq,$$
where $\bm{n}$ is the outward unit normal.

The sparse grid IPDG scheme for \eqref{eq:helmd} developed in \cite{wang2016sparse} is defined as follows. We look for $u_h \in \widehat{\bV}_N^k $, such that
\begin{equation}
\label{sg}
B( u_h, v )= L(v),
\quad  	\forall \,v \in \widehat{\bV}_N^k
\end{equation}
where
\begin{flalign}
\label{eq:bilinear}
B(w, z )= &  \int_{\Omega}  \nabla w \cdot \nabla z\,d\bx
+ \int_{\Omega} cwz\,d\bx -  \sum_{\substack{e \in \Gamma}} \int_{e} \{ \nabla w\} \cdot [ z ]\,ds \notag \\
- & \sum_{\substack{e \in \Gamma}} \int_{e} \{\nabla z\} \cdot [ w ]\,ds
+ \sum_{\substack{e \in \Gamma}}\frac{\sigma}{h} \int_{e} [w] \cdot [z]\,ds,
\end{flalign}
and
\begin{flalign}
\label{eq:rhs}
L(v)= \int_\Omega f v \, d\bx \, ds,
\end{flalign}
where $\sigma$ is a positive penalty parameter, $h=2^{-N}$ is the uniform mesh size in each dimension. Note that the bilinear form $B(\cdot,\cdot)$ is the same as in 
\cite{wheeler1978elliptic, arnold1982interior}, while instead of the expensive full grid piecewise polynomial space $\bV_N^k$, the more efficient sparse grid approximation space $\widehat{\mathbf{V}}_N^k$ is employed, leading to a significant reduction in the size of the linear algebraic system especially when $d$ is large, see \cite{wang2016sparse} for a detailed discussion.

Below, we summarize the theoretical results for  the sparse grid DG methods \eqref{sg}. 
Define the  energy norm of a function $v \in H^2(\Omega_N)$ by
$$ \vertiii{v}^2 :=\sum_{\substack{T \in \Omega_N}}  \int_{T} |\nabla v|^2 \,d\bx \, +\int_{\Omega} c|u|^2\,d\bx +\, \sum_{\substack{e \in \Gamma}}h \int_{e} \left \{  \frac{\partial v}{\partial \bn} \right \}^2\,ds\, + \sum_{\substack{e \in \Gamma}}\frac{1}{h} \int_{e} [v]^2\,ds,$$
where $H^2(\Omega_N)$ denotes the standard broken Sobolev space on grid $\Omega_N$. We also need the following semi-norm.
For any set $L=\left\{i_{1}, \ldots, i_{r}\right\} \subset\{1, \ldots, d\}$, define $L^c$ to be the complement set
of $L$ in $\{1, \ldots, d\}$. For a non-negative integer $\alpha$ and set $L$, we define the mixed derivative semi-norm for a function $v$
$$
|v|_{H^{\alpha, L}\left(\Omega\right)} :=\left\|\left(\frac{\partial^{\alpha}}{\partial x_{i_{1}}^{\alpha}} \cdots \frac{\partial^{\alpha}}{\partial x_{i_{r}}^{\alpha}}\right) v\right\|_{L^{2}\left(\Omega\right)},
$$
and 
$$
|v|_{\mathcal{H}^{q+1}\left(\Omega\right)} :=\max _{1 \leq r \leq d}\left(\max _{L \subset\{1,2, \ldots, d\}}|v|_{H^{q+1, L}\left(\Omega\right)}\right).
$$
The space $\mathcal{H}^{q+1}(\Omega)$ denotes the closure of $C^\infty(\overline{\Omega})$ in the semi-norm $|\cdot|_{\mathcal{H}^{q+1}}$.

We first review several properties related to the projection operators introduced in the previous subsection. We refer the reader to \cite{wang2016sparse,guo2016sparse1} for the proofs.  To avoid unnecessary clutter of constants, the
notation $A\lesssim B$ is used henceforth to represent $A \le C\cdot B$, where the generic constant $C$ is independent of $N$ and the mesh level considered.

\begin{lem}\label{lem:projq}
	Let $\bQ_\bl^k$ be the $L^2$ projection	operator onto the increment space $\bW_\bl^k$. For any $v\in \mathcal{H}^{p+1}(\Omega)$, $k\ge1$,  $1\le q\le\max(p,k)$, $N\ge1$, $d\ge2$, we have
	$$
	\|\bQ_\bl^k (v)\|_{L^2(\Omega)}\lesssim 2^{-(q+1)|\bl|_1}|v|_{\mathcal{H}^{q+1}(\Omega)}.
	$$ 
	
\end{lem}

\begin{lem}[Projection error estimate.]\label{lem:proj}
	Let $\widehat{\bP}_N^k$ be the $L^2$ projection	operator onto the space $\widehat{\bV}_N^k$ introduced in \eqref{eq:sgproj}. For  any $v\in \mathcal{H}^{p+1}(\Omega)$, $k\ge1$, $1\le q\le\max(p,k)$, $N\ge1$, $d\ge2$, we have
	$$
	\vertiii{\widehat{\bP}_N^k(v)-v} \lesssim  N^{d} 2^{-N q}|v|_{\mathcal{H}^{q+1}(\Omega)}.
	$$ 
\end{lem}

The bilinear operator $B(\cdot, \cdot)$ is known to enjoy the following properties.
\begin{lem}[Orthogonality \cite{wang2016sparse}.]
	\label{lem:orth}
	Let $u$ be the exact solution to \eqref{eq:helmd}, and $u_h$ be the numerical solution to \eqref{sg}, then
	$$B( u-u_h, v )=0,   \quad \forall  v \in \widehat{\bV}_N^k.$$
\end{lem}

\begin{lem}[Boundedness and stability \cite{arnold1982interior, Arnold_2002_SIAM_DG}.]
	\label{lem:bound}When $\sigma$ is taken large enough, 
	\begin{align*}B( w, z )& \lesssim \vertiii{w}\cdot \vertiii{z},  && \forall \,w, z \in H^2(\Omega_N);\\
	B( w, w ) &\gtrsim  \vertiii{w}^2,   && \forall \,  v \in \widehat{\bV}_N^k.
	\end{align*}
	
	
\end{lem}

In light of the Lemmas \ref{lem:proj}, \ref{lem:orth} and \ref{lem:bound}, we can prove the following error estimate for the sparse grid IPDG method \eqref{sg}.
\begin{thm}[Convergence \cite{wang2016sparse}.]\label{thm:convergence}
	Let $u$ be the exact solution to \eqref{eq:helmd}, and $u_h$ be the numerical solution to \eqref{sg}. For  $u \in \mathcal{H}^{p+1}(\Omega)$,  $k \geq 1$,  $1 \leq q \leq \min \{p, k\}$,
	$N\geq 1$, $d \geq 2$,  we have
	$$
	\vertiii{u-u_h}   \lesssim  N^d2^{-Nq} |u|_{\mathcal{H}^{q+1}(\Omega)}.
	$$
\end{thm}
This theorem implies a convergence rate of $O(h^{k})$ up to the polylogarithmic term $|\log_2 h|^{d}$ in the energy norm when $u$ is smooth enough. 


Before we proceed, we provide the following estimate, which plays a crucial role in the analysis later. Also note that such an estimate is closely related to the multilevel IPDG method, see e.g. \cite{brix2008multilevel,brix2009multilevel,brix2014multilevel}.
\begin{lem}\label{lem:norm}
	Let $\bP^k_n$ be the $L^2$ projection operator onto the full grid space $\bV_n^k$, $n=0,1\ldots,N$. Then, for any $v\in \bV^k_N$, $k\ge1$, $N\geq 1$, $d \geq 2$, we have
	\begin{equation*}
	\sum_{n=0}^N h_n^{-2}\|(\bP^k_{n}-\bP^k_{n-1})v\|^2_{L^2(\Omega)} \lesssim \vertiii{v}^2,
	\end{equation*}
	where we let $\bP^k_{-1}=0$, and $h_n=2^{-n}$ denote the mesh size of grid $\Omega_n$.
\end{lem}

The proof is given in the Appendix A. By noting that 
$$
\bP_n^k-\bP_{n-1}^k = \sum_{|\bl|_\infty=n} \bQ^k_{\bl},
$$
we further have, for any $v\in \widehat{\mathbf{V}}_N^k$,
\begin{equation}
\label{eq:normeq}
\sum_{n=0}^N 4^n\sum_{\substack{|\bl|_\infty=n\\|\bl|_1\le N}} \|\bQ_{\bl}^k v\|^2_{L^2(\Omega)} \lesssim \vertiii{v}^2.
\end{equation}

\section{Modified IPDG Method on Sparse Grids}
\label{sec:mIPDG}

In this section, we formulate a novel sparse grid IPDG method for solving the Helmholtz equation with variable coefficients. The key idea is to explore the semi-orthogonality property associated with the orthonormal multiwavelet bases.

We start with the following theorem for the constant-coefficient case.
\begin{thm}[Semi-orthogonality.] If the coefficient $c\ge0$ is constant, and two multiwavelet basis functions $v_{\bi,\bl}^{\bj}$ and $v_{\bi',\bl'}^{\bj'}\in\widehat{\mathbf{V}}_N^k$ satisfy $|\max(\bl,\bl')|_1>N$, then 
	\begin{equation}
	\label{eq:so}
	B(v_{\bi,\bl}^{\bj},v_{\bi',\bl'}^{\bj'}) = 0.
	\end{equation}
\end{thm}
\noindent{\it Proof}: Since $v_{\bi,\bl}^{\bj}$, $v_{\bi',\bl'}^{\bj'}\in\widehat{\mathbf{V}}_N^k$, we have 
$$
|\bl|_1\le N,\,|\bl'|_1\le N.
$$
Together with $|\max(\bl,\bl')|_1>N$, we claim that there are two different dimension directions $m_1\ne m_2$ so that $l_{m_1}\ne l'_{m_1}$ and $l_{m_2}\ne l'_{m_2}$. Due to the orthonormal property of the multiwavelet bases, we  follow the argument for the case of the prewavelet bases in \cite{pflaum2016sparse} and prove that \eqref{eq:so} holds.

\hfill $ \blacksquare $

The semi-orthogonality property actually renders a highly desired sparse  structure for the resulting stiffness matrix. On the other hand, when the method \eqref{sg} is applied to \eqref{eq:helmd} with variable coefficients,  semi-orthogonality does not hold anymore, making the  stiffness matrix much denser than that in the constant-coefficient case. To recover the sparse structure  for such a variable-coefficient problem, we propose to modify the bilinear formulation in light of semi-orthogonality, motivated by the work \cite{pflaum1998multilevel,pflaum2016sparse,hartmann2018prewavelet}. In particular, based on $B(\cdot,\cdot)$ introduced in \eqref{eq:bilinear}, we define the following modified bilinear form $B^{so}(\cdot,\cdot):\widehat{\mathbf{V}}_N^k\times \widehat{\mathbf{V}}_N^k \rightarrow \mathbb{R}$
\begin{equation}
\label{eq:def_mbilinear}
B^{so}(\vb,\vbp) = \begin{cases}
B(\vb,\vbp) & |\max(\bl,\bl')|_1\le N,\\
0 & |\max(\bl,\bl')|_1 > N.
\end{cases}
\end{equation}
The sparse grid IPDG weak formulation is modified accordingly as follows. We seek $u_h^{so}\in\sgv$ such that
\begin{equation}
\label{eq:mbilinear}
B^{so}(u_h^{so},v) = L(v),\quad \forall v\in\sgv.
\end{equation}
Note that, by construction, the resulting stiffness matrix by the modified method enjoys the same sparse structure as the constant-coefficient case, leading to computational savings.
Meanwhile, when modifying the bilinear form by means of \eqref{eq:def_mbilinear} we in fact commit a special type of variational crimes. Below, we show that the modified method will generate a numerical solution with the same order accuracy as the original method \eqref{sg} under an extra smoothness assumption of $c$. 

We need the following mixed derivative norm to measure the smoothness for the coefficient $c$. For a function $w$, define the norm
\begin{equation*}
\|w\|_{\mwk(\Omega)}:=\max_{\bzero\le\ba\le\bk+\bone}\left\|  \left(\frac{\partial^{\alpha_1}}{\partial x_1^{\alpha_1}}
\cdots\frac{\partial^{\alpha_d}}{\partial x_d^{\alpha_d}}\right) w\right\|_{L^\infty(\Omega)},
\end{equation*}
and the space $\mwk(\Omega) := \{w\in L^\infty(\Omega): \|w\|_{\mwk(\Omega)}<\infty\}$.
For convenience of illustration, we further introduce several shorthand notation. Whenever given two multi-indexes $\bl$ and $\bl'$, we denote by $\bl^{\max}=\max(\bl,\bl')$, $\bl^{\min}=\min(\bl,\bl')$, and $\Theta_{\bl,\bl'} = \{s\in\{1,2,\ldots,d\}:l_s\ne l_s'\}$.  

We start with the following lemma, which will play a key role in the analysis.
\begin{lem}\label{thm:basis}Let $\vb$ and $\vbp$ be the two multiwavelet basis functions of space $\bW_\bl^k$ and $\bW_{\bl'}^k$, respectively. Denote
	$$\mQ:=\text{supp}(\vb)\cap\text{supp}(\vbp).$$
	Assume $c(\bx)\in\mwk(\Omega)$, then
	\begin{equation}\label{eq:basis_est}\displaystyle
	\left|\int_{\Omega}c(\bx) \vb(\bx)\vbp(\bx) d\bx\right| \lesssim 2^{-\sum_{s\in \Theta_{\bl,\bl'}}l^{\min}_s+(k+1)(l^{\max}_s-l^{\min}_s)}\|c\|_{\mwk(\mQ)} \|\vb\|_{L^2(\mQ)}\|\vbp\|_{L^2(\mQ)}.
	\end{equation}
\end{lem}
\noindent{\it Proof}: If $\mQ$ is empty, then \eqref{eq:basis_est} is trivial. Below, we assume $\mQ$ is nonempty.
Since the multi-dimensional multiwavelet bases introduced in \eqref{eq:basisd} are constructed by tensoring the 1D bases, we can rearrange $\vb(\bx)$ and $\vbp(\bx)$ in the integrand according to $\bl^{\max}$ and $\bl^{\min}$, yielding
\begin{equation}
\int_{\Omega}c(\bx) \vb(\bx)\vbp(\bx) d\bx = \int_{\mQ}c(\bx)v_{\bi^{\max},\bl^{\max}}^{\bj^{\max}}(\bx)
v_{\bi^{\min},\bl^{\min}}^{\bj^{\min}}(\bx)d\bx.\label{eq:equvi}
\end{equation}
%
Assume $\Theta_{\bl,\bl'}=\{m_1,m_2,\ldots,m_r\}$. Note that $l^{\max}_{s}\ge l^{\min}_{s}+1$ iff $s\in\Theta_{\bl,\bl'}$. In
\cite{guo2016sparse1}, we showed that
\begin{equation}
\label{eq:deri}
\left|\int_{\mQ}c(\bx)v_{\bi^{\max},\bl^{\max}}^{\bj^{\max}}(\bx)
v_{\bi^{\min},\bl^{\min}}^{\bj^{\min}}(\bx)d\bx\right| \lesssim  2^{-(k+1)
	\sum_{s\in\Theta_{\bl,\bl'}}l^{\max}_s}\left\|\left(\frac{\partial^{k+1,}}{\partial x_{m_1}^{k+1}}
\cdots\frac{\partial^{k+1}}{\partial x_{m_r}^{k+1}}\right) \left(c\cdot v_{\bi^{\min},\bl^{\min}}^{\bj^{\min}} \right)\right\|_{L^2(\mQ)}.
\end{equation}

Together with the elementary product rule and the fact that  $v_{\bi^{\min},\bl^{\min}}^{\bj^{\min}}$ is a piecewise polynomial of degree $k$ in each dimension, we obtain that
\begin{align*}
\left(\frac{\partial^{k+1,}}{\partial x_{m_1}^{k+1}}
\cdots\frac{\partial^{k+1}}{\partial x_{m_r}^{k+1}}\right) \left(c\cdot v_{\bi^{\min},\bl^{\min}}^{\bj^{\min}} \right)& =(k+1)^r\cdot \left( \frac{\partial}{\partial x_{m_1}}\cdots\frac{\partial}{\partial x_{m_r}}\right)c\cdot
\left( \frac{\partial^{k}}{\partial x_{m_1}^{k}}
\cdots\frac{\partial^{k}}{\partial x_{m_r}^{k}}\right)  v_{\bi^{\min},\bl^{\min}}^{\bj^{\min}} \\
& +\cdots+ \left( \frac{\partial^{k+1}}{\partial x_{m_1}^{k+1}}
\cdots\frac{\partial^{k+1}}{\partial x_{m_r}^{k+1}}\right)c\cdot v_{\bi^{\min},\bl^{\min}}^{\bj^{\min}}.
\end{align*}
A direct application of the inverse inequality leads to 
\begin{equation}\label{eq:deriv2}
\left\|\left(\frac{\partial^{k+1,}}{\partial x_{m_1}^{k+1}}
\cdots\frac{\partial^{k+1}}{\partial x_{m_r}^{k+1}}\right) \left(c\cdot v_{\bi^{\min},\bl^{\min}}^{\bj^{\min}} \right)\right\|_{L^2(\mQ)} \lesssim 2^{k
	\sum_{s\in\Theta_{\bl,\bl'}}l^{\min}_s} \|c\|_{
	\mwk(\mQ)}\|v_{\bi^{\min},\bl^{\min}}^{\bj^{\min}}\|_{L^2(\mQ)}.
\end{equation}
Also note that $\|v_{\bi^{\max},\bl^{\max}}^{\bj^{\max}}\|_{L^2(\mQ)}=1$, and by definition, $$\|v_{\bi^{\max},\bl^{\max}}^{\bj^{\max}}\|_{L^2(\mQ)}\|v_{\bi^{\min},\bl^{\min}}^{\bj^{\min}}\|_{L^2(\mQ)} = \|\vb\|_{L^2(\mQ)}\|\vbp\|_{L^2(\mQ)}.$$
This, together with  \eqref{eq:equvi}, \eqref{eq:deri} and \eqref{eq:deriv2}, completes the proof.

\hfill $ \blacksquare $

\begin{lem}\label{thm:space}Assume $c(\bx)\in\mwk(\Omega)$. If $w(\bx)\in \bW^k_{\bl}$, and $z(\bx)\in \bW^k_{\bl'}$, then
	$$
	\left|\int_{\Omega} c(\bx)w(\bx)z(\bx)d\bx\right| \lesssim 2^{-\sum_{s\in \Theta_{\bl,\bl'}}l^{\min}_s+(k+1)(l^{\max}_s-l^{\min}_s)}\|c\|_{\mwk(\Omega)} \|w\|_{L^2(\Omega)}\|z\|_{L^2(\Omega)}.
	$$
\end{lem}
\noindent{\it Proof}: Note that, for 
$w(\bx)\in\bW^k_{\bl}$ and $z(\bx)\in\bW^k_{\bl'}$, we have
\begin{align*}
w(\bx)& = \sum_{\bj\in \Theta_{\bl}}\sum_{\bone\le\bi\le\bk+\bone} w_{\bi,\bl}^\bj v_{\bi,\bl}^\bj(\bx),\\
z(\bx)& = \sum_{\bj\in \Theta_{\bl'}}\sum_{\bone\le\bi\le\bk+\bone} z_{\bi,\bl'}^\bj v_{\bi,\bl'}^\bj(\bx).
\end{align*}
Then, based on Lemma \ref{thm:basis}, the proof immediately follows the same procedure as in Lemma 4.6 in \cite{pflaum2016sparse} and hence is omitted for brevity.

\hfill $ \blacksquare $

The above lemma directly leads to the following estimate, which is useful when proving the boundedness and stability of the modified bilinear form $B^{so}(\cdot,\cdot)$.
\begin{lem}
	\label{lem:wz}
	Assume $c(\bx)\in\mwk(\Omega)$, $k\ge1$.
	If $w(\bx)\in \bW^k_{\bl}$, and $z(\bx)\in \bW^k_{\bl'}$, then
	$$
	\left|\int_{\Omega} c(\bx)w(\bx)z(\bx)d\bx\right| \lesssim  2^{-\sum_{s\in \Theta_{\bl,\bl'}}l^{\max}_s}\|c\|_{\mwk(\Omega)} \|w\|_{L^2(\Omega)}\|z\|_{L^2(\Omega)}.
	$$
\end{lem}
\noindent{\it Proof}: Since $k\ge1$,
$$
-\sum_{s\in \Theta_{\bl,\bl'}}l^{\min}_s+(k+1)(l^{\max}_s-l^{\min}_s)\le -\sum_{s\in \Theta_{\bl,\bl'}}l^{\min}_s+(l^{\max}_s-l^{\min}_s) = - \sum_{s\in \Theta_{\bl,\bl'}}l^{\max}_s.
$$

\hfill $ \blacksquare $


Now, we are ready to establish the boundedness and stability of the modified bilinear form $B^{so}(\cdot,\cdot)$.
\begin{thm}[Boundedness and stability with semi-orthogonality.]\label{thm:bound} Assume $c(\bx)\in\mwk(\Omega)$, $k\ge1$, $d\ge2$. There exists 
	an integer $N_0$, such that
	\begin{align*}
	B^{so}(w,z)& \lesssim \vertiii{w}\cdot\vertiii{z},\\
	B^{so}(w,w)&  \gtrsim \vertiii{w},
	\end{align*}	
	where $w,z\in \widehat{\mathbf{V}}_N^k$ with $N\ge N_0$.
\end{thm}

\noindent{\it Proof}: We first show that, for  $w,z\in \widehat{\mathbf{V}}_N^k$,
\begin{equation}
\label{eq:diff}
|B(w,z)-B^{so}(w,z)|\lesssim 2^{-N_d}N^d\vertiii{w}\cdot\vertiii{z},
\end{equation}
where
\begin{equation*}
N_d:=\begin{cases}
N & d\le 4,\\
\frac{2}{d-2}N & d>4. 
\end{cases}
\end{equation*}
Note that 
\begin{align}
\begin{split}\label{eq:bilinear_diff}
&|B(w,z)-B^{so}(w,z)| \\
&= \left|\sum_{\substack{ |\bl|_1\le N, |\bl'|_1\le N\\ |\bl^{\max}|_1>N}}\int_{\Omega}c(\bx)\bQ^k_{\bl}(w)\bQ^k_{\bl'}(z)d\bx\right|\\
&\le \sum_{\substack{ |\bl|_1\le N, |\bl'|_1\le N\\ |\bl^{\max}|_1>N}}\left|\int_{\Omega}c(\bx)\bQ^k_{\bl}(w)\bQ^k_{\bl'}(z)d\bx\right|\\
&\lesssim \sum_{\substack{ |\bl|_1\le N, |\bl'|_1\le N\\ |\bl^{\max}|_1>N}}
2^{-\sum_{s\in \Theta_{\bl,\bl'}}l^{\max}_s}\|c\|_{\mwk(\Omega)} \|_{L^2(\Omega)}\|\bQ^k_{\bl}(w)\|_{L^2(\Omega)}\|\bQ^k_{\bl'}(z)\|_{L^2(\Omega)}\\
& =\|c\|_{\mwk(\Omega)} \sum_{\substack{ |\bl|_1\le N, |\bl'|_1\le N\\ |\bl^{\max}|_1>N}} 2^{-|\bl|_\infty-|\bl'|_{\infty}-\sum_{s\in \Theta_{\bl,\bl'}}l^{\max}_s}2^{|\bl|_\infty}\|\bQ^k_{\bl}(w)\|_{L^2(\Omega)}2^{|\bl'|_\infty}\|\bQ^k_{\bl'}(z)\|_{L^2(\Omega)}
\end{split}
\end{align}
where we have used Lemma \ref{lem:wz}.
The rest of the proof follows the same procedure as in  Theorem 5.4 in \cite{pflaum2016sparse}. The only difference is that we need the Lemma \ref{thm:space} and estimate \eqref{eq:normeq} to account for the discontinuous approximation space. 
%
For the completeness of the paper, we choose to provide the proof.

First, in \cite{pflaum2016sparse}, it was shown that if $|\bl|_1\le N$, $|\bl'|_1\le N$, and $|\bl^{\max}|_1>N$, then
$$
|\bl|_\infty+|\bl'|_{\infty}+\sum_{s\in \Theta_{\bl,\bl'}}l_s^{\max}\ge |\bl^{\max}|_1 - N + N_d.
$$
Then,
\begin{align}
\begin{split}\label{eq:bilinear_diff2}
&\sum_{\substack{ |\bl|_1\le N, |\bl'|_1\le N\\ |\bl^{\max}|_1>N}} 2^{-|\bl|_\infty-|\bl'|_{\infty}-\sum_{s\in \Theta_{\bl,\bl'}}l^{\max}_s}2^{|\bl|_\infty}\|\bQ^k_{\bl}(w)\|_{L^2(\Omega)}2^{|\bl'|_\infty}\|\bQ^k_{\bl'}(z)\|_{L^2(\Omega)}\\
&\le 2^{-N_d}\sum_{\substack{ |\bl|_1\le N, |\bl'|_1\le N\\ |\bl^{\max}|_1>N}}2^{-(|\bl^{\max}|_1-N)}2^{|\bl|_\infty}\|\bQ^k_{\bl}(w)\|_{L^2(\Omega)}2^{|\bl'|_\infty}\|\bQ^k_{\bl'}(z)\|_{L^2(\Omega)}\\
& = 2^{-N_d}\sum_{\gamma=1}^N2^{-\gamma}\sum_{\substack{ |\bl|_1\le N, |\bl'|_1\le N\\ |\bl^{\max}|_1=N+\gamma}}2^{|\bl|_\infty}\|\bQ^k_{\bl}(w)\|_{L^2(\Omega)}2^{|\bl'|_\infty}\|\bQ^k_{\bl'}(z)\|_{L^2(\Omega)}\\
&\le 2^{-N_d} \sqrt{\sum_{\gamma=1}^N2^{-\gamma}\sum_{\substack{ |\bl|_1\le N, |\bl'|_1\le N\\ |\bl^{\max}|_1=N+\gamma}}4^{|\bl|_\infty}\|\bQ^k_{\bl}(w)\|^2_{L^2(\Omega)}}\sqrt{\sum_{\gamma=1}^N2^{-\gamma}\sum_{\substack{ |\bl|_1\le N, |\bl'|_1\le N\\ |\bl^{\max}|_1=N+\gamma}}4^{|\bl'|_\infty}\|\bQ^k_{\bl'}(z)\|^2_{L^2(\Omega)}},
\end{split}
\end{align}
where we have used the Cauchy-Schwarz inequality.
Note that
\begin{align}
\label{eq:norm_sum}
\sum_{\gamma=1}^N2^{-\gamma}\sum_{\substack{ |\bl|_1\le N, |\bl'|_1\le N\\ |\bl^{\max}|_1=N+\gamma}}4^{|\bl|_\infty}\|\bQ^k_{\bl}(w)\|^2_{L^2(\Omega)}= \sum_{n=0}^N \sum_{\substack{ |\bl|_1\le N\\ |\bl|_\infty=n}}4^n\|\bQ^k_{\bl}(w)\|^2_{L^2(\Omega)}\sum_{\gamma=1}^N 2^{-\gamma}\sum_{\substack{|\bl'|_1\le N\\ |\bl^{\max}|_1=N+\gamma}}1.
\end{align}
We bound the right-hand side of above equation as follows. First, we need the inequality proved in \cite{pflaum2016sparse} that, for any $\bl$ with $|\bl|_1\le N$, 
$$
\sum_{\gamma=1}^N 2^{-\gamma}\sum_{\substack{|\bl'|_1\le N\\ |\bl^{\max}|_1=N+\gamma}}1 \le C N^d,
$$
where $C$ is a constant independent of $N$ and $\bl$. Together with the estimate \eqref{eq:normeq}, we derive that 
\begin{equation}
\label{eq:bilinear_diff3}
\sum_{n=0}^N \sum_{\substack{ |\bl|_1\le N\\ |\bl|_\infty=n}}4^n\|\bQ^k_{\bl}(w)\|^2_{L^2(\Omega)}\sum_{\gamma=1}^N 2^{-\gamma}\sum_{\substack{|\bl'|_1\le N\\ |\bl^{\max}|_1=N+\gamma}}1
\lesssim N^d \vertiii{w}.
\end{equation}
Combining \eqref{eq:bilinear_diff}, \eqref{eq:bilinear_diff2}, \eqref{eq:norm_sum} and \eqref{eq:bilinear_diff3}, we prove \eqref{eq:diff}. Due to the boundedness and stability of $B(\cdot,\cdot)$, we complete the proof of the theorem.  

\hfill $ \blacksquare $

Now, we are ready to establish the main result of the paper. 

\begin{thm}[Convergence with semi-orthogonality.] \label{thm:main}
	Let $u$ be the exact solution to the Helmholtz equation \eqref{eq:helmd}, and $u^{so}_h$ be the numerical solution to the modified IPDG formulation with semi-orthogonality \eqref{eq:mbilinear}. For  $u \in \mathcal{H}^{p+1}(\Omega)$, $k \geq 1$,  $c\in\mwk(\Omega)$, $q=\max(p,k)$, $d \geq 2$, 
	$N\geq N_0$, where $N_0$ is given in Theorem~\ref{thm:bound},  we have
	$$
	\vertiii{u-u_h^{so}}    \lesssim  N^d 2^{-qN} |u|_{\mathcal{H}^{q+1}(\Omega)}.
	$$
\end{thm}
\noindent{\it Proof}: Following the theory of variational crime, see \cite{grossmann1992numerik,brenner2007mathematical}, we consider the decomposition of the error 
\begin{equation}
\label{eq:err_decom}
\vertiii{u-u^{so}_h} \le \vertiii{u-\widehat{\bP}^k_N(u)} + \vertiii{u^{so}_h-\widehat{\bP}^k_N(u)}.
\end{equation}
For the second term on the right-hand side, due to stability of $B^{so}(\cdot,\cdot)$ in $\widehat{\mathbf{V}}_N^k$ from Theorem \ref{thm:bound}, we have
\begin{align}
\label{eq:variational_crime}
\begin{split}
\vertiii{u^{so}_h-\widehat{\bP}^k_N{u}}& \lesssim \sup_{w\in \widehat{\mathbf{V}}_N^k}\frac{B^{so}(u^{so}_h-\widehat{\bP}^k_n(u),w)}{\vertiii{w}}\\
&= \sup_{w\in \widehat{\mathbf{V}}_N^k}\frac{B^{so}(u^{so}_h,w)-B^{so}(\widehat{\bP}^k_n(u),w)}{\vertiii{w}}\\
& = \sup_{w\in \widehat{\mathbf{V}}_N^k}\frac{B(u_h,w)-B^{so}(\widehat{\bP}^k_n(u),w)}{\vertiii{w}}\\
& = \sup_{w\in \widehat{\mathbf{V}}_N^k}\frac{B(u,w)-B^{so}(\widehat{\bP}^k_n(u),w)}{\vertiii{w}}\\
& = \sup_{w\in \widehat{\mathbf{V}}_N^k}\frac{B(u -\widehat{\bP}^k_n(u) ,w) +B(\widehat{\bP}^k_n(u) ,w) -B^{so}(\widehat{\bP}^k_n(u),w)}{\vertiii{w}}\\
& \lesssim \vertiii{u-\widehat{\bP}^k_N} + \sup_{w\in \widehat{\mathbf{V}}_N^k}\frac{B(\widehat{\bP}^k_n(u) ,w) -B^{so}(\widehat{\bP}^k_n(u),w)}{\vertiii{w}},
\end{split}
\end{align}
where $u_h$ denotes the numerical solution by the original IPDG method  \eqref{eq:bilinear}. In the derivation, we have also used orthogonality  and boundedness of $B(\cdot,\cdot)$. Note that the first term on the right-hand side is the projection error and has been estimated in Lemma \ref{lem:proj}; while the second term  measures the effect of modifying the bilinear form. For the rest of the proof, we show that
\begin{equation}
\label{eq:diff_err2}
|B(\widehat{\bP}^k_N(u),v)-B^{so}(\widehat{\bP}^k_N(u),v)|\lesssim N^{d}2^{-(q+1)N}|u|_{\mathcal{H}^{q+1}(\Omega)}\vertiii{v}.
\end{equation} 
Since $\bQ_\bl^k\left(\widehat{\bP}^k_N(u)\right) = \bQ_\bl^k(u)$ for $|\bl|_1\le N$, we derive that
\begin{align}
\begin{split}
\label{eq:diff_err}
&|B(\widehat{\bP}^k_N(u),v)-B^{so}(\widehat{\bP}^k_N(u),v)|\\
&\le \sum_{\substack{ |\bl|_1\le N, |\bl'|_1\le N\\ |\bl^{\max}|_1>N}}\left|\int_{\Omega}c(\bx)\bQ^k_{\bl}(u)\bQ^k_{\bl'}(v)d\bx\right|\\
&\lesssim \sum_{\substack{ |\bl|_1\le N, |\bl'|_1\le N\\ |\bl^{\max}|_1>N}}2^{-\sum_{s\in \Theta_{\bl,\bl'}}l^{\min}_s+(k+1)(l^{\max}_s-l^{\min}_s)}\|c\|_{\mwk(\Omega)} \|\bQ_\bl^k (u)\|_{L^2(\Omega)}\|\bQ_{\bl'}^k(v)\|_{L^2(\Omega)}\\
&\lesssim \|c\|_{\mwk(\Omega)}|u|_{\mathcal{H}^{q+1}(\Omega)}\sum_{\substack{ |\bl|_1\le N, |\bl'|_1\le N\\ |\bl^{\max}|_1>N}}2^{-\sum_{s\in \Theta_{\bl,\bl'}}l^{\min}_s+(k+1)(l^{\max}_s-l^{\min}_s)} 2^{-(q+1)|\bl|_1}\|\bQ_{\bl'}^k(v)\|_{L^2(\Omega)},
\end{split}
\end{align}
where we have used Lemmas \ref{thm:space} and \ref{lem:projq}.
Notice that
\begin{align*}\sum_{s\in \Theta_{\bl,\bl'}}l^{\min}_s+(k+1)(l^{\max}_s-l^{\min}_s)+(q+1)|\bl|_1&\ge \sum_{s\in \Theta_{\bl,\bl'}}(k+1)(l^{\max}_s-l^{\min}_s)+(q+1)|\bl|_1 \\
&\ge (q+1)|\bl^{\max}|_1.
\end{align*}
Then, we have
\begin{align*}
|B(\widehat{\bP}^k_N(u),v)-B^{so}(\widehat{\bP}^k_N(u),v)|
&\lesssim \|c\|_{\mwk(\Omega)}|u|_{\mathcal{H}^{q+1}(\Omega)}\sum_{\substack{ |\bl|_1\le N, |\bl'|_1\le N\\ |\bl^{\max}|_1>N}}2^{-(q+1)|\bl^{\max}|_1} \|\bQ_{\bl'}^k(v)\|_{L^2(\Omega)}\\
&=2^{-(q+1)N}\|c\|_{\mwk(\Omega)}|u|_{\mathcal{H}^{q+1}(\Omega)}\sum_{\gamma=1}^{N}2^{-(q+1)\gamma}\sum_{\substack{ |\bl|_1\le N, |\bl'|_1\le N\\ |\bl^{\max}|_1=N+\gamma}} \|\bQ_{\bl'}^k(v)\|_{L^2(\Omega)}.
\end{align*}
Similar to \eqref{eq:bilinear_diff}-\eqref{eq:bilinear_diff3}, we make use of the Cauchy-Schwarz inequality and  bound the sum on the right-hand side as
\begin{align*}
\sum_{\gamma=1}^{N}2^{-(q+1)\gamma}\sum_{\substack{ |\bl|_1\le N, |\bl'|_1\le N\\ |\bl^{\max}|_1=N+\gamma}} \|\bQ_{\bl'}^k(v)\|_{L^2(\Omega)}
\le\sum_{\gamma=1}^{N}2^{-\gamma}\sum_{\substack{ |\bl|_1\le N, |\bl'|_1\le N\\ |\bl^{\max}|_1=N+\gamma}} \|\bQ_{\bl'}^k(v)\|_{L^2(\Omega)}
\lesssim N^{d}\vertiii{v}.
\end{align*}
The proof of \eqref{eq:diff_err2} is complete.

By combining \eqref{eq:err_decom}, \eqref{eq:variational_crime}, \eqref{eq:diff_err2}, and Lemma \ref{lem:proj} about the projection error estimate, we complete the proof.

\hfill $ \blacksquare $

\begin{rem}
	$\sup_{w\in \widehat{\mathbf{V}}_N^k}\frac{B(\widehat{\bP}^k_Nu ,w) -B^{so}(\widehat{\bP}^k_Nu,w)}{\vertiii{w}}$ quantifies the variational crime from modifying the bilinear form and is indeed one order higher than the projection error. From the numerical results presented in the next section, we will see that if the coefficient $c$ and the solution $u$ are sufficiently smooth, then the modified sparse grid IPDG method will generate almost the same numerical results as the original method, while the resulting linear system is much sparser, leading to additional computational savings when solving \eqref{eq:helmd} with variable coefficients.
\end{rem}

\begin{rem} The proposed framework can be extended to other variable-coefficient problems if a similar estimate as in Lemma \ref{thm:basis} is derived. For instance, for
	Poisson's equation with variable coefficients $-\nabla\cdot(\bK(\bx)\nabla u )=f$, under some condition of $\bK$, it is possible to estimate
	$$
	\left|\int_{\Omega}\nabla \vb(\bx)\cdot\left(\bK(\bx)\nabla\vbp(\bx)\right) d\bx\right|
	$$
	and the boundary terms in the  IPDG formulation with $|\bl|_1\le N$, $|\bl|_1\le N$ and $|\bl^{\max}|_1>N$,	
	and devise a modified sparse grid method with semi-orthogonality in an attempt to sparsify the stiffness matrix and save cost. We leave the investigation to the future work.
\end{rem}

\input{numerical}

\section{Conclusion}
\label{sec:Con}
In this paper, we developed a modified sparse grid IPDG method using semi-orthogonality for solving the Helmholtz equation with variable coefficients. 
The original IPDG features a sparse finite element space which scales as $O(h^{-1}|\log_2 h|^{d-1})$ for $d$-dimensional problems,  translating into a significant cost reduction when $d$ is large. On the other hand, when applied to the variable-coefficient problem, the method suffers a dense stiffness matrix, impeding its efficiency advantage to some extent. Based on the semi-orthogonality property associated with the orthonormal multiwavelet bases, the IPDG bilinear form is modified aiming to sparsify the resulting linear algebraic matrix. A numerical analysis demonstrates that the error incurred by the modification is one order higher than the projection error, making it negligible in the simulations. Numerical results  in up to six dimensions are shown to validate the analysis and  demonstrate that the modified method enjoys sparser stiffness matrix than the original one, leading to extra computational savings.
%
Future work includes the study of adaptive algorithm for less regular solutions as well as extension of the method to  other types of high-dimensional variable-coefficient equations.

\begin{appendices}
	
	\section{Proof of Lemma \ref{lem:norm}}
	In the proof, we need the averaging projection operator $\mathcal{A}$ proposed in \cite{brix2008multilevel,brix2009multilevel}, which  decomposes the space $\bV_N^k$ into the conforming subspace $\bV_N^{c}$ and nonconforming part $\bV_N^{nc}$, i.e,
	$$\bV_N^{c} := \mathcal{A}\bV_N^k \subset\bV_N^k\cap H^1_0(\Omega),\quad \bV_N^{nc} := \bV_N^k - \bV_N^{c}.$$
	Furthermore, $\mathcal{A}$ satisfies the Jackson estimate
	\begin{equation}
	\label{ap:jackson}
	\vertiii{\mathcal{A}v}_a\lesssim \vertiii{v}_a,\quad \vertiii{(\mathcal{I}-\mathcal{A})v}_a \lesssim h\vertiii{v}_a,\quad v\in \bV_N^k,
	\end{equation}
	where the energy norm $\vertiii{v}_a$ is defined by  
	$$
	\vertiii{v}_a :=\sum_{\substack{T \in \Omega_N}}  \int_{T} |\nabla v|^2 \,d\bx \,  + \sum_{\substack{e \in \Gamma}}\frac{1}{h} \int_{e} [v]^2\,ds,\quad v\in \bV_{N}^k,
	$$
	and  $\mathcal{I}$ denotes the identity operator. The construction and detailed analysis of $\mathcal{A}$ were established  in \cite{brix2008multilevel,brix2009multilevel,brix2014multilevel}. 
	Clearly, $\vertiii{v}_a\le\vertiii{v}$.
	
	For any $v\in \bV_N^k$, we denote by $w := \mathcal{A}(v)$ and $z: = (\mathcal{I}-\mathcal{A})v$. We first prove that 
	\begin{equation}
	\label{ap:w}
	\sum_{n=0}^{N} h_n^{-2}\|(\bP^k_n-\bP^k_{n-1})w\|^2_{L^2(\Omega)}\lesssim |w|_{H^1(\Omega)}.
	\end{equation}
	Note that there is a standard nested relation for conforming finite element spaces:
	$$
	\bV^{c}_0\subset \bV^{c}_1\subset\cdots \bV^{c}_N\subset H_0^1(\Omega),
	$$
	where the space $\bV^{c}_n$ defined on mesh $\Omega_n$ is a subspace of $\bV^{k}_n$, $n=0,\ldots,N$.
	Following the idea of the BPX multilevel precondtioner \cite{bramble1990parallel}, we define the Ritz projection operator $\bR_n: \bV^{c}_N \rightarrow \bV^{c}_n$
	$$
	a(\bR_n(w),r) = a(w,r) \quad \forall r\in  \bV^{c}_n,\quad n = 1,\ldots,N,
	$$
	where 
	$$a(w,r) = \int_\Omega \nabla w\cdot \nabla r d\bx,$$
	and we let $\bR_{-1}=0$.
	Denote $w_n:=(\bR_n-\bR_{n-1})w\in 	\bV_n^c$. Then, $w=\sum_{n=0}^N w_n$. Since $\bV_n^{c}\subset \bV_n^k$,
	\begin{equation}
	\label{ap:orth}
	(\bP_l^k - \bP_{l-1}^k) w_n=0,\quad l> n.
	\end{equation}
	In addition, the Ritz operator enjoys the property
	$$\bR_{n-1}(\bR_n-\bR_{n-1}) = 0.$$ 
	Together with the standard estimate of $\bR_{n-1}$ by the conforming finite element analysis, 
	we have
	\begin{equation}
	\label{ap:wn}
	\|w_n\|_{L^2(\Omega)}= \|(\mathcal{I}-\bR_{n-1})w_n\|_{L^2(\Omega)}\lesssim h_n |w_n|_{H^1(\Omega)}.
	\end{equation}
	Then, we derive 
	\begin{align*}
	\sum_{l=0}^{N} h_l^{-2}\|(\bP^k_l-\bP^k_{l-1})w\|^2_{L^2(\Omega)}& =\sum_{l=0}^{N} h_l^{-2}\|(\bP^k_l-\bP^k_{l-1})\sum_{n=0}^{N}w_n\|^2_{L^2(\Omega)}\\
	&=\sum_{l=0}^{N} h_l^{-2}\|\sum_{n=0}^{N}(\bP^k_l-\bP^k_{l-1})w_n\|^2_{L^2(\Omega)}\\
	& = \sum_{l=0}^{N} h_l^{-2}\sum_{n,m=0}^{N}((\bP^k_l-\bP^k_{l-1})w_n,(\bP^k_l-\bP^k_{l-1})w_m)\\
	& = \sum_{l=0}^{N} h_l^{-2}\sum_{n,m=l}^{N}((\bP^k_l-\bP^k_{l-1})w_n,(\bP^k_l-\bP^k_{l-1})w_m) \tag{ \text{due to \eqref{ap:orth}}} \\
	& = \sum_{n,m=0}^{N}\sum_{l=0}^{\min(n,m)} h_l^{-2}((\bP^k_l-\bP^k_{l-1})w_n,(\bP^k_l-\bP^k_{l-1})w_m)\\
	& \le \sum_{n,m=0}^{N}\sum_{l=0}^{\min(n,m)} h_l^{-2}\|(\bP^k_l-\bP^k_{l-1})w_n\|_{L^2(\Omega)}\|(\bP^k_l-\bP^k_{l-1})w_m\|_{L^2(\Omega)}\\
	& \lesssim \sum_{n,m=0}^{N}\sum_{l=0}^{\min(n,m)} h_l^{-2}\|w_n\|_{L^2(\Omega)}\|w_m\|_{L^2(\Omega)}\\
	&\lesssim \sum_{n,m=0}^{N}\sum_{l=0}^{\min(n,m)} h_l^{-2}h_nh_m|w_n|_{H^1(\Omega)}|w_m|_{H^1(\Omega)}\tag{ \text{due to \eqref{ap:wn}}}\\
	&\lesssim \sum_{n,m=0}^{N} \left(\frac{1}{2}\right)^{|n-m|}|w_n|_{H^1(\Omega)}|w_m|_{H^1(\Omega)}.
	\end{align*}
	By an elementary linear algebra result, see \cite{brenner2007mathematical}, together with the identity $a(w_n,w_m)=0$, $n\ne m$, we finally get
	$$
	\sum_{n,m=0}^{N} \left(\frac{1}{2}\right)^{|n-m|}|w_n|_{H^1(\Omega)}|w_m|_{H^1(\Omega)}\lesssim \sum_{n=0}^N|w_n|_{H^1(\Omega)}^2 = |w|_{H^1(\Omega)}^2,
	$$
	and \eqref{ap:w} is proved. Together with \eqref{ap:w} and the Jackson estimate \eqref{ap:jackson}, we immediately obtain
	\begin{equation}
	\label{ap:w1}
	\sum_{n=0}^{N} h_n^{-2}\|(\bP^k_n-\bP^k_{n-1})w\|^2_{L^2(\Omega)}\lesssim\vertiii{w}_a\lesssim \vertiii{v}_a.
	\end{equation}
	
	Now we derive the estimate for $z$. A direct application of the Poincar\'e-Friedrichs type inequality for space $\mathbf{V}_N^k$, see \cite{arnold1982interior}, yields
	$$
	\|z\|_{L^2(\Omega)} \lesssim \vertiii{z}_a.
	$$
	This, together with the Jackson estimate \eqref{ap:jackson}, gives
	\begin{align}
	\label{ap:z}
	\sum_{n=0}^{N} h_n^{-2}\|(\bP^k_n-\bP^k_{n-1})z\|^2_{L^2(\Omega)}\lesssim  \sum_{n=0}^{N} h_n^{-2}\|z\|^2_{L^2(\Omega)}\lesssim  \sum_{n=0}^{N} h_n^{-2}\vertiii{z}_a^2\lesssim \left(h^2\sum_{n=0}^{N} h_n^{-2}\right)\vertiii{v}_a^2\lesssim \vertiii{v}_a^2.
	\end{align}
	Combining the estimates \eqref{ap:w1} and \eqref{ap:z}, we complete the proof as
	\begin{align*}
	\sum_{n=0}^{N} h_n^{-2}\|(\bP^k_n-\bP^k_{n-1})v\|^2_{L^2(\Omega)}&\lesssim \sum_{n=0}^{N} h_n^{-2}\|(\bP^k_n-\bP^k_{n-1})w\|^2_{L^2(\Omega)} +\sum_{n=0}^{N} h_n^{-2}\|(\bP^k_n-\bP^k_{n-1})z\|^2_{L^2(\Omega)}\\&\lesssim \vertiii{v}_a^2\le \vertiii{v}^2.
	\end{align*}
	
	\hfill $\blacksquare$
\end{appendices}

\bibliographystyle{abbrv}
\bibliography{ref_cheng,ref_cheng_2,ref_guo}

\end{document}

%% file: numerical.tex
\section{Numerical Results}
\label{sec: SGDG numerical}
In this section, we provide numerical results to demonstrate the performance of the modified sparse grid IPDG method as well as verify the analysis established for simulating the Helmholtz equations with variable coefficients. The penalty parameter $\sigma$ is set to be empirical values  $\sigma = 5\cdot k\cdot d$.

\begin{figure}[htp]
\begin{center}
\subfigure[]{\includegraphics[width=.47\textwidth]{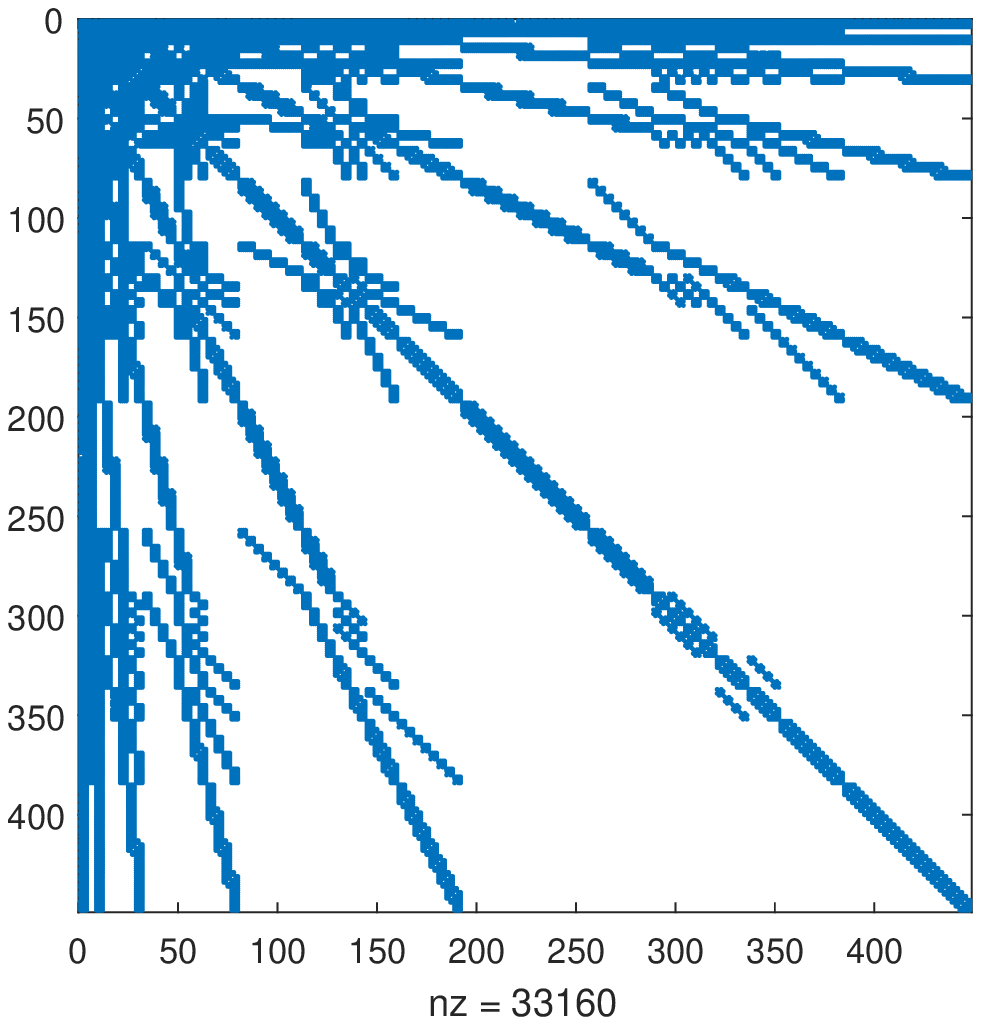}}
\subfigure[]{\includegraphics[width=.47\textwidth]{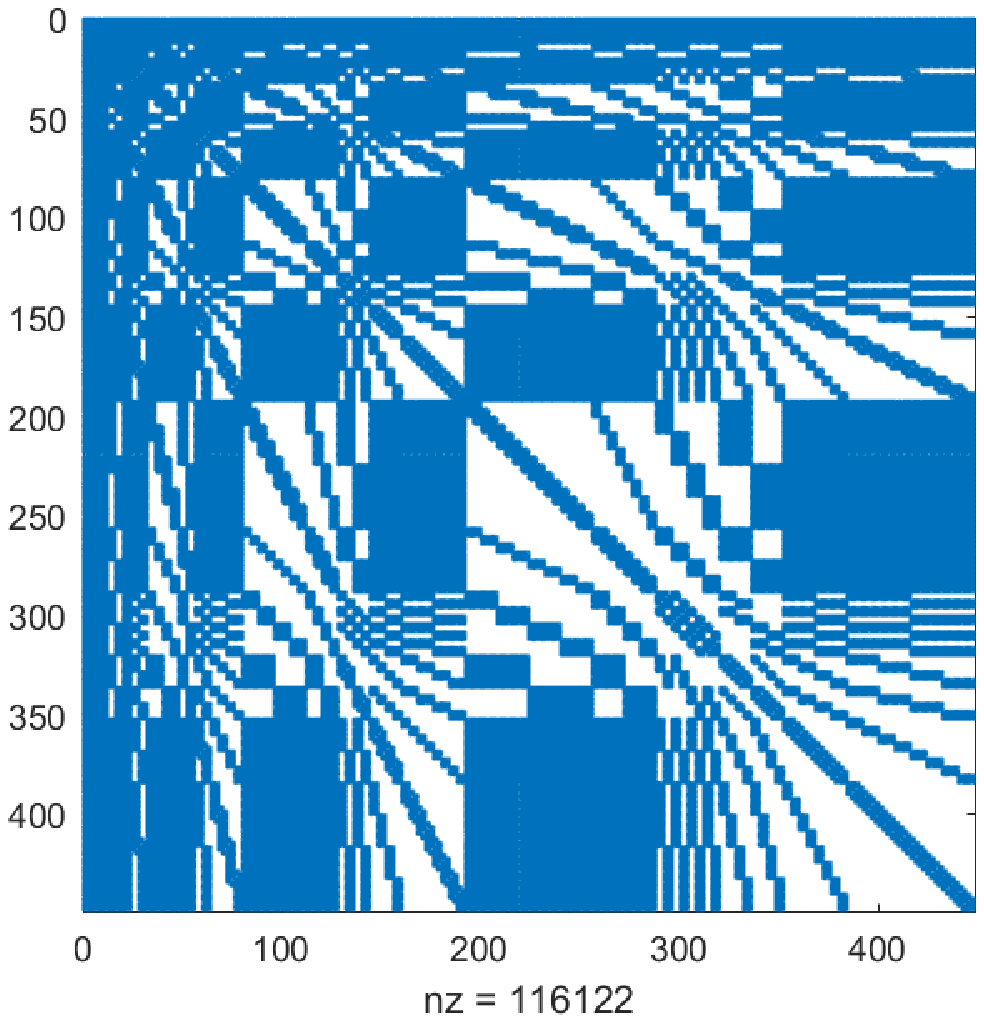}}
\end{center}
  \caption{Example \ref{exa:1}.  The sparsity patterns of the stiffness matrices by  (a) the modified sparse grid IPDG method and (b) the original sparse grid IPDG method. $N=5$, $k=1$, $d=2$.}
 \label{fig:spy2d}
\end{figure}

\begin{exa}\label{exa:1} We solve the Helmholtz equation \eqref{eq:helmd} with the smooth coefficient
$$c(\bx) = \prod_{m=1}^d (1-x_m^2),$$
up to $d=6$. The right-hand $f$ is chosen such that 
$$u(\bx)=\prod_{m=1}^d \sin(\pi x_m).$$

\end{exa}
This problem was solved in \cite{hartmann2018prewavelet} by a continuous Ritz-Galerkin discretization on sparse grids using prewavelets in conjunction with semi-orthogonality.	
We first consider the case of $d=2$ and compare and contrast the performance of both modified and original sparse grid IPDG methods in terms of accuracy and efficiency. In Table \ref{table:exa12d}, we summarize the convergence study for both methods with $k=1,\,2,\,3$. It is observed that, for this problem, the $L^2$ and $H^1$ errors given by the two methods are virtually the same, and the associated convergence rates are close to $k+1$-th order for the $L^2$ error and $k$-th order for the $H^1$ error. Such an observation
validates the convergence analysis in Theorem \ref{thm:main} that the variational crime incurred by the proposed modification of the bilinear form is one order higher than the projection error and hence becoming negligible in the simulation when the smoothness requirements are fulfilled. 
To demonstrate the efficiency advantage of the modified method, we provide the numbers of nonzero entries of the stiffness matrices as well as the conditional numbers for both methods with various mesh levels $N$ and polynomial degrees $k$ in Table \ref{table:exa12d_nonzeros}. Further, in Figure \ref{fig:spy2d}, we plot the sparsity patterns of the matrices with $N=5$ and $k=1$. One can see that  the modified method enjoys much sparser stiffness matrices than the original one, while the condition numbers remain almost unchanged. More specifically, we compute the order of sparsity defined by $O_s = \log(\text{NNZ})/\log(\text{DOF})$, where NNZ is referred to as  
the number of nonzero elements. Note that NNZ scales as $O(\text{DOF}^{1.5})$ for the modified method; while it scales as $O(\text{DOF}^{1.7})$ for the original method. Such reduction leads to great computational savings when assembling as well as solving the algebraic linear system.

\begin{table}
\caption{Numerical errors and orders of accuracy for Example \ref{exa:1} computed by the modified and the original IPDG methods with $k=1,\, 2,\,3$. $d=2$. 
}
\vspace{2 mm}
\centering
\begin{tabular}{ |c | c c c c |c c  c c |}
\hline
 & \multicolumn{4}{c|}{Modified IPDG method }  & \multicolumn{4}{c|}{Original IPDG method }\\\hline

$N$& $L^2$ error & order & $H^1$ error & order& $L^2$ error & order& $H^1$ error & order\\\hline
&\multicolumn{8}{c|}{$k=1$ } 		\\\hline
2	&	6.12E-02	&		&	6.88E-01	&		&	6.12E-02	&		&	6.88E-01	&		\\\hline
3	&	1.62E-02	&	1.92	&	3.37E-01	&	1.03	&	1.62E-02	&	1.92	&	3.37E-01	&	1.03	\\\hline
4	&	4.02E-03	&	2.01	&	1.66E-01	&	1.02	&	4.02E-03	&	2.01	&	1.66E-01	&	1.02	\\\hline
5	&	9.99E-04	&	2.01	&	8.26E-02	&	1.01	&	9.99E-04	&	2.01	&	8.26E-02	&	1.01	\\\hline
6	&	2.52E-04	&	1.99	&	4.11E-02	&	1.01	&	2.52E-04	&	1.99	&	4.11E-02	&	1.01	\\\hline
&\multicolumn{8}{c|}{$k=2$ } 			\\\hline													
2	&	1.52E-03	&		&	5.27E-02	&		&	1.52E-03	&		&	5.27E-02	&		\\\hline
3	&	2.27E-04	&	2.74	&	1.34E-02	&	1.98	&	2.27E-04	&	2.74	&	1.34E-02	&	1.98	\\\hline
4	&	3.51E-05	&	2.70	&	3.35E-03	&	2.00	&	3.51E-05	&	2.70	&	3.35E-03	&	2.00	\\\hline
5	&	5.27E-06	&	2.74	&	8.37E-04	&	2.00	&	5.27E-06	&	2.74	&	8.37E-04	&	2.00	\\\hline
6	&	7.61E-07	&	2.79	&	2.09E-04	&	2.00	&	7.61E-07	&	2.79	&	2.09E-04	&	2.00	\\\hline
&\multicolumn{8}{c|}{$k=3$ } 	\\\hline
2	&	9.17E-05	&		&	3.53E-03	&		&	9.17E-05	&		&	3.53E-03	&		\\\hline
3	&	6.03E-06	&	3.93	&	4.37E-04	&	3.01	&	6.03E-06	&	3.93	&	4.37E-04	&	3.01	\\\hline
4	&	3.84E-07	&	3.97	&	5.43E-05	&	3.01	&	3.84E-07	&	3.97	&	5.43E-05	&	3.01	\\\hline
5	&	2.43E-08	&	3.98	&	6.78E-06	&	3.00	&	2.43E-08	&	3.98	&	6.78E-06	&	3.00	\\\hline
6	&	1.54E-09	&	3.98	&	8.47E-07	&	3.00	&	1.54E-09	&	3.98	&	8.47E-07	&	3.00	\\\hline

\end{tabular}
\label{table:exa12d}
\end{table}

\begin{table}
	\caption{Sparsity and condition number of the stiffness matrices for the modified and original IPDG methods. $k =1,\,2,\,3$. $d=2$. DOF is the number of degrees of freedom used for the sparse grid IPDG methods. NNZ is the number of nonzero elements in the stiffness matrix. $O_s$=log(NNZ)/ log(DOF). $\kappa_2$ is the condition number. 
	}
	\vspace{2 mm}
	\centering
	\begin{tabular}{ |c  c | c c c |c c  c  |}
		\hline
		 \multicolumn{2}{|c|}{$d=2$ } & \multicolumn{3}{c|}{Modified IPDG method } & \multicolumn{3}{c|}{Original IPDG method } \\\hline
		$N$ & DOF & NNZ & $O_s$ & $\kappa_2$ &  NNZ & $O_s$&  $\kappa_2$\\\hline
		& &\multicolumn{6}{c|}{$k=1$ } 		\\\hline
2	&	32	&	592	&	1.84	&	8.23E+01	&	976	&	1.99	&	8.23E+01	\\\hline
3	&	80	&	2608	&	1.80	&	3.49E+02	&	5424	&	1.96	&	3.49E+02	\\\hline
4	&	448	&	9808	&	1.51	&	1.40E+03	&	26192	&	1.67	&	1.40E+03	\\\hline
5	&	1024	&	33160	&	1.50	&	5.54E+03	&	116122	&	1.68	&	5.54E+03	\\\hline
6	&	2304	&	103968	&	1.49	&	2.20E+04	&	493154	&	1.69	&	2.20E+04	\\\hline
	& &\multicolumn{6}{c|}{$k=2$ } 		\\\hline														
2	&	180	&	2976	&	1.54	&	3.51E+02	&	4920	&	1.64	&	3.51E+02	\\\hline
3	&	432	&	12864	&	1.56	&	1.37E+03	&	27120	&	1.68	&	1.37E+03	\\\hline
4	&	1008	&	48132	&	1.56	&	5.35E+03	&	131076	&	1.70	&	5.35E+03	\\\hline
5	&	2304	&	162324	&	1.55	&	2.11E+04	&	580352	&	1.71	&	2.11E+04	\\\hline
6	&	5184	&	509616	&	1.54	&	8.37E+04	&	2460726	&	1.72	&	8.37E+04	\\\hline
	& &\multicolumn{6}{c|}{$k=3$ }		\\\hline														

2	&	128	&	9216	&	1.88	&	8.53E+02	&	15616	&	1.99	&	8.53E+02	\\\hline
3	&	768	&	39952	&	1.59	&	3.29E+03	&	85760	&	1.71	&	3.29E+03	\\\hline
4	&	1792	&	149264	&	1.59	&	1.28E+04	&	413952	&	1.73	&	1.28E+04	\\\hline
5	&	4096	&	505072	&	1.58	&	5.04E+04	&	1848064	&	1.73	&	5.04E+04	\\\hline
6	&	9216	&	1587200	&	1.56	&	2.00E+05	&	7869696	&	1.74	&	2.00E+05	\\\hline

	\end{tabular}
	\label{table:exa12d_nonzeros}
\end{table}

\begin{figure}[htp]
	\begin{center}
		\subfigure[]{\includegraphics[width=.47\textwidth]{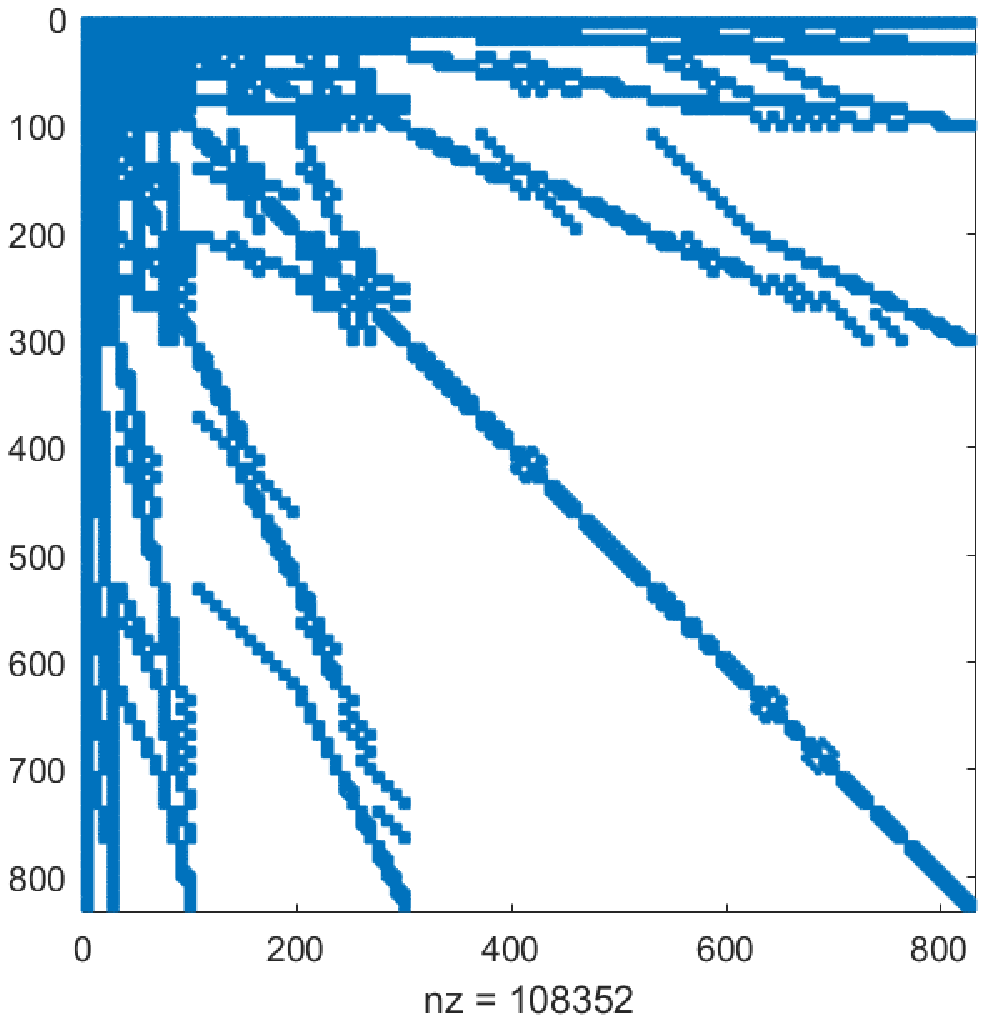}}
		\subfigure[]{\includegraphics[width=.47\textwidth]{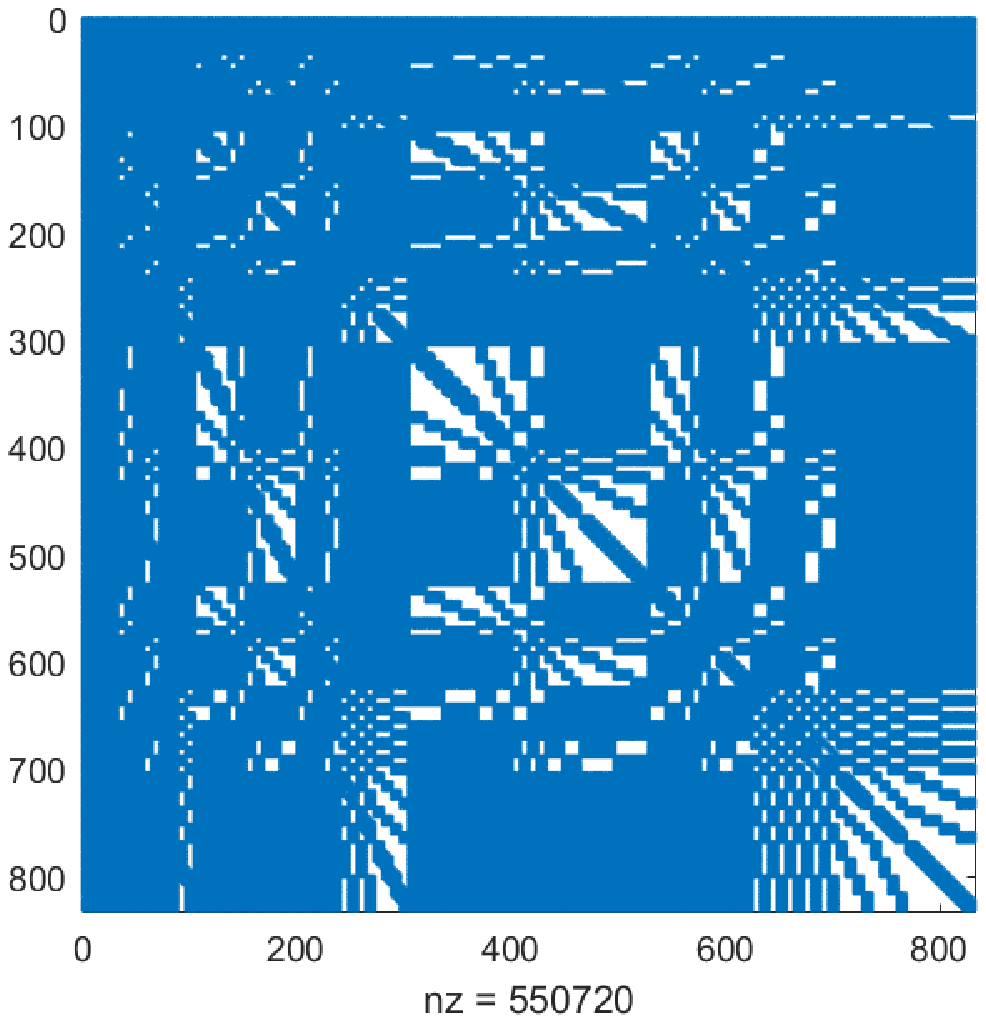}}
	\end{center}
	\caption{Example \ref{exa:1}.  The sparsity patterns of the stiffness matrices by  (a) the modified sparse grid IPDG method and (b) the original sparse grid IPDG method. $N=4$, $k=1$, $d=3$. }
	\label{fig:spy3d}
\end{figure}

We then consider the case of $d=3$ and perform a similar comparison between the two IPDG methods. We provide the convergence study results in Table \ref{table:exa13d} with $k=1,\,2,\,3,$ including the $L^2$ and $H^1$ errors and the associated orders of accuracy for both methods. As in the above $d=2$ case, both methods give  almost the same numerical errors. $k+1$-th and $k$-th order of accuracy is observed for the $L^2$ and $H^1$ errors, respectively, which agrees with the error estimates established in Theorem \ref{thm:main}. In Figure \ref{fig:spy3d}, we plot the sparsity patterns of the stiffness matrices by two methods with $N=4$ and $k=1$. As expected, the stiffness matrix by the modified method is much sparser than that by the original method. Furthermore, we provide DOF, NNZ and the associated order of sparsity $O_s$ for both methods in Table \ref{table:exa13d_nonzeros}. It is observed that NNZ scales as $O(\text{DOF}^{1.7})$ for the modified method; while it scales as $O(\text{DOF}^{1.95})$ for the original method. 

\begin{table}
	\caption{Numerical errors and orders of accuracy for Example \ref{exa:1} computed by the modified and the original IPDG methods with $k=1,\, 2,\,3$. $d=3$. 
	}
	\vspace{2 mm}
	\centering
	\begin{tabular}{ |c | c c c c |c c  c c |}
		\hline
		& \multicolumn{4}{c|}{Modified IPDG method }  & \multicolumn{4}{c|}{Original IPDG method }\\\hline
		
		$N$& $L^2$ error & order & $H^1$ error & order& $L^2$ error & order& $H^1$ error & order\\\hline
		&\multicolumn{8}{c|}{$k=1$ } 		\\\hline
		2	&	1.30E-01	&		&	9.45E-01	&		&	1.30E-01	&		&	9.45E-01	&		\\\hline
		3	&	3.54E-02	&	1.87	&	4.41E-01	&	1.10	&	3.54E-02	&	1.87	&	4.41E-01	&	1.10	\\\hline
		4	&	8.88E-03	&	1.99	&	2.07E-01	&	1.09	&	8.88E-03	&	1.99	&	2.07E-01	&	1.09	\\\hline
		5	&	2.16E-03	&	2.04	&	9.90E-02	&	1.07	&	2.16E-03	&	2.04	&	9.90E-02	&	1.07	\\\hline
		6	&	5.42E-04	&	1.99	&	4.82E-02	&	1.04	&	5.42E-04	&	1.99	&	4.82E-02	&	1.04	\\\hline
		&\multicolumn{8}{c|}{$k=2$ } 			\\\hline													
2	&	1.13E-03	&		&	4.54E-02	&		&	1.13E-03	&		&	4.54E-02	&		\\\hline	
3	&	2.10E-04	&	2.43	&	1.19E-02	&	1.94	&	2.10E-04	&	2.43	&	1.19E-02	&	1.94	\\\hline	
4	&	3.73E-05	&	2.49	&	3.01E-03	&	1.98	&	3.73E-05	&	2.49	&	3.01E-03	&	1.98	\\\hline	
5	&	6.15E-06	&	2.60	&	7.54E-04	&	2.00	&	6.15E-06	&	2.60	&	7.54E-04	&	2.00	\\\hline	
6	&	9.71E-07	&	2.66	&	1.88E-04	&	2.00	&	9.71E-07	&	2.66	&	1.88E-04	&	2.00	\\\hline

		&\multicolumn{8}{c|}{$k=3$ } 	\\\hline
		2	&	1.12E-04	&		&	1.30E-03	&		&	1.12E-04	&		&	1.30E-03	&		\\\hline
		3	&	7.40E-06	&	3.92	&	1.39E-04	&	3.23	&	7.40E-06	&	3.92	&	1.39E-04	&	3.23	\\\hline
		4	&	4.72E-07	&	3.97	&	1.58E-05	&	3.14	&	4.72E-07	&	3.97	&	1.58E-05	&	3.14	\\\hline
		5	&	2.99E-08	&	3.98	&	1.86E-06	&	3.08	&	2.99E-08	&	3.98	&	1.86E-06	&	3.08	\\\hline
		6	&	1.89E-09	&	3.98	&	2.26E-07	&	3.05	&	1.89E-09	&	3.98	&	2.26E-07	&	3.05	\\\hline

	\end{tabular}
	\label{table:exa13d}
\end{table}

\begin{table}
\caption{Sparsity of the stiffness matrices for the modified and original IPDG methods. $k =1,\,2,\,3$. $d=3$. DOF is the number of degrees of freedom used for the sparse grid IPDG methods. NNZ is the number of nonzero elements in the stiffness matrix. $O_s$=log(NNZ)/ log(DOF). 
}
\vspace{2 mm}
\centering
\begin{tabular}{ |c  c | c c  |c c    |}
	\hline
	\multicolumn{2}{|c|}{$d=3$ }  & \multicolumn{2}{c|}{Modified IPDG method } & \multicolumn{2}{c|}{Original IPDG method } \\\hline
	$N$ & DOF & NNZ & $O_s$  &  NNZ & $O_s$ \\\hline
	& &\multicolumn{4}{c|}{$k=1$ } 		\\\hline
2	&	104	&	4.31E+03	&	1.80	&		1.05E+04	&	1.99		\\\hline	
3	&	304	&	2.34E+04	&	1.76	&		8.26E+04	&	1.98		\\\hline	
4	&	832	&	1.08E+05	&	1.72	&		5.51E+05	&	1.97		\\\hline	
5	&	2176	&	4.47E+05	&	1.69	&		3.27E+06	&	1.95		\\\hline	
6	&	5504	&	1.69E+06	&	1.67	&		1.80E+07	&	1.94		\\\hline

	& &\multicolumn{4}{c|}{$k=2$ } 		\\\hline														
2	&	351	&	4.93E+04	&	1.84	&	1.19E+05	&		1.99		\\\hline	
3	&	1026	&	2.65E+05	&	1.80	&	9.38E+05	&		1.98		\\\hline	
4	&	2808	&	1.22E+06	&	1.76	&	6.26E+06	&		1.97		\\\hline	
5	&	7344	&	5.03E+06	&	1.73	&	3.73E+07	&		1.96		\\\hline	
6	&	18576	&	1.91E+07	&	1.71	&	2.06E+08	&		1.95		\\\hline

	& &\multicolumn{4}{c|}{$k=3$ }		\\\hline														
	
2	&	832	&	2.76E+05	&	1.86	&		6.69E+05	&	1.99		\\\hline	
3	&	2432	&	1.48E+06	&	1.82	&		5.26E+06	&	1.99		\\\hline	
4	&	6656	&	6.81E+06	&	1.79	&		3.51E+07	&	1.97		\\\hline	
5	&	17408	&	2.81E+07	&	1.76	&		2.10E+08	&	1.96		\\\hline	
6	&	44032	&	1.07E+08	&	1.73	&		1.16E+09	&	1.95		\\\hline

\end{tabular}
\label{table:exa13d_nonzeros}
\end{table}

Last, we summarize the simulation results for $d=4,\,5,\,6$ in Tables \ref{table:exa14d}-\ref{table:exa16d}, respectively. In particular, we report the $L^2$ errors, the associated orders of accuracy, NNZ, and the orders of sparsity $O_s$ for the modified method only. The observation is similar to that in the cases of $d=2,\,3$. The modified method is able to achieve the expected orders of convergence, while slight order reduction is observed due to the polylogarithmic term appearing in the error estimate.
Note that unlike the modified method, the original method suffers the almost fully populated stiffness matrix, making its simulation exceed our computation resource limit. Currently, we are allowed to handle  matrices with NNZ no greater than  $1E+10$. Based on the comparison drawn in $d=2,\,3$ together with the Theorems \ref{thm:convergence} and \ref{thm:main}, we believe the errors and the associated convergence rates of the two methods are also comparable. 

\begin{table}
	\caption{The $L^2$ errors, orders of accuracy and sparsity of the stiffness matrices for the modified IPDG methods. $k =1,\,2,\,3$. $d=4$. DOF is the number of degrees of freedom used for the sparse grid IPDG methods. NNZ is the number of nonzero elements in the stiffness matrix. $O_s$=log(NNZ)/ log(DOF). 
	}
	\vspace{2 mm}
	\centering
	\begin{tabular}{ |c  c | c c  |c c    |}
		\hline
			\multicolumn{2}{|c|}{$d=4$ }  & \multicolumn{4}{c|}{Modified IPDG method }  \\\hline
		$N$ & DOF & $L^2$ error & order  &  NNZ & $O_s$ \\\hline
		\multicolumn{6}{|c|}{$k=1$ } 		\\\hline

2	&	304	&	1.49E-01	&		&	2.76E+04	&	1.79	\\\hline	
3	&	1008	&	7.35E-02	&	1.02	&	1.78E+05	&	1.75	\\\hline	
4	&	3072	&	2.15E-02	&	1.77	&	9.71E+05	&	1.72	\\\hline	
5	&	8832	&	5.62E-03	&	1.93	&	4.69E+06	&	1.69	\\\hline	
6	&	24320	&	1.44E-03	&	1.97	&	2.06E+07	&	1.67	\\\hline	
\multicolumn{6}{|c|}{$k=2$ } 		\\\hline
2	&	1539	&	1.28E-03	&		&	7.04E+05	&	1.83	\\\hline	
3	&	5103	&	2.25E-04	&	2.51	&	4.52E+06	&	1.80	\\\hline	
4	&	15552	&	3.94E-05	&	2.51	&	2.47E+07	&	1.76	\\\hline	
5	&	44712	&	6.63E-06	&	2.57	&	1.19E+08	&	1.74	\\\hline	
6	&	123120	&	1.09E-06	&	2.61	&	5.25E+08	&	1.71	\\\hline	
\multicolumn{6}{|c|}{$k=3$ } 		\\\hline

2	&	4864	&	7.83E-05	&		&	7.02E+06	&	1.86	\\\hline	
3	&	16128	&	5.30E-06	&	3.88	&	4.51E+07	&	1.82	\\\hline	
4	&	49152	&	3.44E-07	&	3.95	&	2.46E+08	&	1.79	\\\hline

	\end{tabular}
	\label{table:exa14d}
\end{table}

\begin{table}
	\caption{The $L^2$ errors, orders of accuracy and sparsity of the stiffness matrices for the modified IPDG methods. $k =1,\,2,\,3$. $d=5$. DOF is the number of degrees of freedom used for the sparse grid IPDG methods. NNZ is the number of nonzero elements in the stiffness matrix. $O_s$=log(NNZ)/ log(DOF). 
	}
	\vspace{2 mm}
	\centering
	\begin{tabular}{ |c  c | c c  |c c    |}
		\hline
					\multicolumn{2}{|c|}{$d=5$ }  & \multicolumn{4}{c|}{Modified IPDG method }  \\\hline
		$N$ & DOF & $L^2$ error & order  &  NNZ & $O_s$ \\\hline
		\multicolumn{6}{|c|}{$k=1$ } 		\\\hline
		
2	&	832	&	1.30E-01	&		&	1.60E+05	&	1.78	\\\hline
3	&	3072	&	9.37E-02	&	0.47	&	1.20E+06	&	1.74	\\\hline
4	&	10272	&	4.10E-02	&	1.19	&	7.52E+06	&	1.71	\\\hline
5	&	32064	&	1.31E-02	&	1.65	&	4.15E+07	&	1.69	\\\hline
6	&	95104	&	3.71E-03	&	1.81	&	2.07E+08	&	1.67	\\\hline
	\multicolumn{6}{|c|}{$k=2$ } 		\\\hline

2	&	6318	&	1.06E-03	&		&	9.22E+06	&	1.83	\\\hline
3	&	23328	&	1.97E-04	&	2.42	&	6.89E+07	&	1.79	\\\hline
4	&	78003	&	3.60E-05	&	2.45	&	4.32E+08	&	1.77	\\\hline
5	&	243486	&	6.29E-06	&	2.52	&	2.38E+09	&	1.74	\\\hline
6	&	722196	&	1.08E-06	&	2.54	&	1.19E+10	&	1.72	\\\hline
	\multicolumn{6}{|c|}{$k=3$ } 		\\\hline
2	&	26624	&	6.79E-05	&		&	1.64E+08	&	1.86	\\\hline
3	&	98304	&	4.69E-06	&	3.85	&	1.22E+09	&	1.82	\\\hline
4	&	328704	&	3.10E-07	&	3.92	&7.67E+09	&	1.79	\\\hline

	\end{tabular}
	\label{table:exa15d}
\end{table}

\begin{table}
	\caption{The $L^2$ errors, orders of accuracy and sparsity of the stiffness matrices for the modified IPDG methods. $k =1,\,2$. $d=6$. DOF is the number of degrees of freedom used for the sparse grid IPDG methods. NNZ is the number of nonzero elements in the stiffness matrix. $O_s$=log(NNZ)/ log(DOF). 
	}
	\vspace{2 mm}
	\centering
	\begin{tabular}{ |c  c | c c  |c c    |}
		\hline
					\multicolumn{2}{|c|}{$d=6$ }  & \multicolumn{4}{c|}{Modified IPDG method }  \\\hline
		$N$ & DOF & $L^2$ error & order  &  NNZ & $O_s$ \\\hline
		\multicolumn{6}{|c|}{$k=1$ } 		\\\hline
		
2	&	2176	&	1.03E-01	&		&	8.78E+05	&	1.78	\\\hline	
3	&	8832	&	8.71E-02	&	0.24	&	7.48E+06	&	1.74	\\\hline	
4	&	32064	&	5.70E-02	&	0.61	&	5.29E+07	&	1.71	\\\hline	
5	&	107712	&	2.33E-02	&	1.29	&	3.27E+08	&	1.69	\\\hline	
6	&	341504	&	7.85E-03	&	1.57	&	1.81E+09	&	1.67	\\\hline	

	\multicolumn{6}{|c|}{$k=2$ } 		\\\hline
2	&	24786	&	8.47E-04	&		&	1.14E+08	&	1.83	\\\hline
3	&	100602	&	1.67E-04	&	2.35	&	9.69E+08	&	1.80	\\\hline
4	&	365229	&	3.13E-05	&	2.41	&	6.85E+09	&	1.77	\\\hline
		
	\end{tabular}
	\label{table:exa16d}
\end{table}
